\pdfoutput=1
\newif\ifpreprint%
\preprinttrue   

\ifpreprint
  \documentclass[12pt]{article}
\else
  \RequirePackage{fix-cm}
  \documentclass[smallextended]{svjour3}
\fi

\usepackage{hyperref} 

\ifpreprint
\hypersetup{colorlinks,
            citecolor=mdarkgreen,
            linkcolor=mmred}

\usepackage{fullpage} 
\usepackage{parskip}

\usepackage[small,labelfont={bf,sf}]{caption}
\usepackage{etoolbox}
\patchcmd{\thebibliography}{\advance\leftmargin\labelsep}
  {\labelsep=0.7cm \advance\leftmargin\labelsep}{}{}

\usepackage{abstract}

\usepackage{multirow}
\AtBeginDocument{}%
\AtBeginDocument{}%
\AtBeginDocument{}%

\usepackage{amsthm}

\makeatletter
\def\thm@space@setup{%
  \thm@preskip=5mm plus 2mm minus 2mm
  \thm@postskip=\thm@preskip 
}
\makeatother

\newtheorem{theorem}{Theorem}

\theoremstyle{definition}
\newtheorem{definition}{Definition}[section]

\theoremstyle{remark}

\else
    \hypersetup{colorlinks,
            citecolor=black,
            linkcolor=blue}
\fi

\usepackage[utf8]{inputenc}
\usepackage[OT1]{fontenc}

\usepackage{amsmath,amssymb,amsfonts}
\usepackage{enumerate}

\usepackage[english]{babel}
\usepackage{graphicx}
\usepackage[caption=false]{subfig}
 \usepackage{times}

 \usepackage{placeins} 
\usepackage{mathtools} %
\usepackage[ruled,vlined,linesnumbered]{algorithm2e}

\usepackage{siunitx}
\usepackage{xargs}                      
\usepackage[table]{xcolor}  
 \usepackage{csquotes} 
\usepackage{csvsimple} 
\usepackage{booktabs}
\usepackage[para,online,flushleft]{threeparttable}

\usepackage{tikz}
\usepackage{pgfplots}
\usepackage{pgfplotstable}
\usepackage{gincltex} 
\usepackage{pgfmath,pgffor}
\pgfplotsset{compat=newest} 
\pgfplotsset{
      table/search path={benchmark_data},
    }
\usetikzlibrary{shapes.geometric,calc,arrows,intersections, external, matrix, fit, quotes, shapes}

\tikzset{
    sbox/.style={
        draw,
        rectangle split,
        rectangle split parts=2,
        text centered,
    },
    sbox+/.style={
        sbox,
        rectangle split every empty part={},
        rectangle split empty part width=\widthof{#1},
        rectangle split empty part height=\heightof{#1},
        rectangle split empty part depth=\depthof{#1},
    },
}

\ifpreprint
  \newcommand{\eqdef}{\coloneqq}
  \newcommand{\reqdef}{\eqqcolon}
  \newcommand{\pkg}[1]{\texttt{#1}}
\else
  \newcommand{\eqdef}{=}
  \newcommand{\reqdef}{=}
  \newcommand{\pkg}[1]{#1} 
\fi

\newcommand{\tx}{\tilde{x}}
\newcommand{\ts}{\tilde{s}}

\newcommand{\hx}{\hat{x}}
\newcommand{\hs}{\hat{s}}

\newcommand{\norm}[1]{\left\lVert#1\right\rVert}
\newcommand{\abs}[1]{\left|#1\right|}

\newcommand{\mc}[1]{\mathcal{#1}}
\newcommand{\mb}[1]{\mathbb{#1}}

\DeclareMathOperator*{\argmin}{argmin}

\makeatletter
\newcommand{\removelatexerror}{\let\@latex@error\@gobble}
\makeatother

\definecolor{mmred}{HTML}{e0024c}
\definecolor{mred}{rgb}{0.863,0.129,0.302}
\definecolor{mgreen}{rgb}{0,0.549,0}   
\definecolor{mblue}{rgb}{0,0.392,0.871} 
\definecolor{mpurple}{rgb}{0.4,0.0,0.4} 
\definecolor{morange}{rgb}{1,0.4,0} 
\definecolor{mteal}{rgb}{0.216,0.784,0.671} 
\definecolor{mdarkgreen}{RGB}{11, 91, 35}
\definecolor{oxfordblue}{RGB}{2,46,95}

\newcommand{\myCellColor}{red!25}

\newcommand{\OptProblem}[5][]{
  \begin{aligned}\label{#1}
    \begin{array}{ll}
    \underset{#3}{\rm{#2}}&\hspace{-0ex}#4\vspace{.5ex}\\
    \rm{subject~to}&#5
  \end{array}
\end{aligned}
}
\newcommand{\OptProblemExtra}[5][]{
  \begin{aligned}\label{#1}
    \begin{array}{lll}
    \underset{#3}{\rm{#2}}&\hspace{-0ex}#4\vspace{.5ex}\\
    \rm{subject~to}&#5
  \end{array}
\end{aligned}
}

\DeclareMathOperator{\prox}{prox}

\newcommand{\svec}[1]{\mathrm{svec}({#1})}

\newcommand{\inv}{^{-1}}
\newcommand{\project}[1]{{\Pi}_{#1}}

\renewcommand{\Re}{\mathbb{R}}
\newcommand{\Sym}[1][n]{\mathbb{S}^{#1}}
\newcommand{\Psd}[1][n]{\mathbb{S}_{+}^{#1}}
\newcommand{\clique}[1]{\mathcal{C}_{#1}}
\newcommand{\sep}[1]{\mathrm{sep}({#1})}

\newcommand{\eg}{{\it e.g.}\xspace}

\definecolor{mosekcolor}{rgb}{0.863,0.129,0.302}
\definecolor{scscolor}{rgb}{0,0.549,0}   
\definecolor{cosmocolor}{rgb}{0,0.392,0.871} 
\definecolor{cosmocdcolor}{rgb}{1,0.4,0} 
\definecolor{osqpcolor}{rgb}{0.4,0.0,0.4} 
\definecolor{gurobicolor}{RGB}{11, 91, 35}

\newcommand{\gmeanRatios}{0.701}
\newcommand{\minRatio}{0.407} 
\newcommand{\minRatioName}{mcp500-2}

\ifpreprint

  \title{COSMO: A conic operator splitting method \\for convex conic problems}

  \author{ Michael Garstka$^*$
  \and Mark Cannon$^*$
  \and Paul Goulart
  \thanks{The authors are with the Department of Engineering Science, University of Oxford, Oxford, OX1 3PJ, UK. Email: {\{michael.garstka, mark.cannon, paul.goulart\}@eng.ox.ac.uk}}
  }
\else
  \title{COSMO: A conic operator splitting method for convex conic problems\thanks{Preliminary versions of this work appeared in~\cite{garstka_2019,Garstka_2020}. MG is supported by the Clarendon Scholarship.}}

  \author{ Michael Garstka \and
  Mark Cannon \and
  Paul Goulart}

  \institute{
  Michael Garstka \at
  University of Oxford, Parks Road, Oxford, OX1 3PJ, UK\\
  \email{michael.garstka@eng.ox.ac.uk}           
  \and
  Mark Cannon \at
  University of Oxford, Parks Road, Oxford, OX1 3PJ, UK\\
  \email{mark.cannon@eng.ox.ac.uk}           
  \and
  Paul Goulart \at
  University of Oxford, Parks Road, Oxford, OX1 3PJ, UK\\
  \email{paul.goulart@eng.ox.ac.uk}           
  }

  \date{Received: date / Accepted: date}

 \journalname{Mathematical Programming Computation}
\fi

\begin{document}

\maketitle


\begin{abstract}
This paper describes the Conic Operator Splitting Method (COSMO) solver, an operator splitting algorithm for convex optimisation problems with quadratic objective function and conic constraints.
At each step the algorithm alternates between solving a quasi-definite linear system with a constant coefficient matrix and a projection onto convex sets.
The low per-iteration computational cost makes the method particularly efficient for large problems, e.g.\ semidefinite programs that arise in portfolio optimisation, graph theory, and robust control.
Moreover, the solver uses chordal decomposition techniques and a new clique merging algorithm to effectively exploit sparsity in large, structured semidefinite programs.
A number of benchmarks against other state-of-the-art solvers for a variety of problems show the effectiveness of our approach.
Our Julia implementation is open-source, designed to be extended and customised by the user, and is integrated into the Julia optimisation ecosystem.

\ifpreprint \else
\fi

\end{abstract}

\section{Introduction}
We consider convex optimisation problems in the form
\begin{equation}\label{eqn:BasicProblem}
   \begin{array}{lll}
    \mbox{minimize}   & f(x) &\\
    \mbox{subject to} & g_i(x)  \leq 0, & i=1,\ldots,l \\
                      & h_i(x) = 0, & i=1,\ldots,k,
  \end{array}
\end{equation}
where we assume that both the objective function $f : \Re^n \to \Re$ and the inequality constraint functions $g_i : \Re^n \to \Re$ are convex, and that  the equality constraints $h_i(x) := a_i^\top  x -b_i$ are affine.
We will denote an optimal solution to this problem (if it exists) as $x^*$.
Convex optimisation problems feature heavily in a wide range of research areas and industries, including problems in machine learning~\cite{Cortes_1995}, finance~\cite{Boyd_2017a}, optimal control~\cite{Boyd_1994}, and operations research~\cite{Borm_2001}.
Concrete examples of problems fitting the general form~\eqref{eqn:BasicProblem} include linear programming (LP), convex quadratic programming (QP), second-order cone programming (SOCP), and semidefinite programming (SDP) problems.
Methods to solve each of these standard problem classes are well known  and a number of open- and closed-source solvers are widely available.
However, the trend for data and training sets of increasing size in decision making problems and machine learning poses a challenge for state-of-the-art software.

Algorithms for LPs were first used to solve military planning and allocation problems in the 1940s~\cite{Dantzig_1963}. In 1947 Danzig developed the simplex method that solves LPs by searching for the optimal solution along the vertices of the inequality polytope.
Extensions to the method led to the general field of \emph{active-set methods}~\cite{Wolfe_1959} that are able to solve both LPs and QPs, and
which search for an optimal point by iteratively constructing a set of active constraints.
Although often efficient in practice, a major theoretical drawback is that the worst-case complexity increases exponentially with the problem size~\cite{Wright_1999}.

The most common approach taken by modern convex solvers is the \emph{
interior-point method}~\cite{Wright_1999}, which stems from Karmarkar's original projective algorithm~\cite{Karmarkar_1984}, and is able to solve both LPs and QPs in polynomial time.
Interior point methods have since been extended to problems with positive semidefinite (PSD) constraints in~\cite{Helmberg_1996} and~\cite{Alizadeh_1998}.
The primal-dual interior point methods apply variants of Newton's method to iteratively find a solution to a set of optimality KKT conditions.
At each iteration the algorithm alternates between a Newton step that involves factoring a Jacobian matrix and a line search to determine the magnitude of the step to ensure a feasible iterate.
Most notably, the Mehrotra predictor-corrector method in~\cite{Mehrotra_1992} forms the basis of several implementations because of its strong practical performance~\cite{Wright_1997}.
However, interior-point methods typically do not scale well for very large problems since the Jacobian matrix has to be calculated and factored at each step.

Two main approaches to overcome this limitation are active research areas.
Firstly, a renewed focus on first-order methods with computationally cheaper per-iteration-cost, and secondly the exploitation of sparsity in the problem data.
First-order methods are known to handle larger problems well at the expense of reduced accuracy compared to interior-point methods.
In the 1960s Everett~\cite{Everett_1963} proposed a dual decomposition method that allows one to decompose a separable objective function, making each iteration cheaper to compute.
Augmented Lagrangian methods by Miele~(\cite{Miele_1971,Miele_1971a,Miele_1972}), Hestenes~\cite{Hestenes_1969}, and Powell~\cite{Powell_1969} are more robust and helped to remove the strict convexity conditions on problems, while losing the decomposition property.
By splitting the objective function, the alternating direction method of multipliers (ADMM), first described in~\cite{Gabay_1975}, allowed the advantages of dual decomposition to be combined with the superior convergence and robustness of augmented Lagrangian methods.
Subsequently, it was shown that ADMM can be analysed from the perspective of monotone operators and that it is a special case of Douglas-Rachford splitting~\cite{Eckstein_1989} as well as of the proximal point algorithm in~\cite{Rockafellar_1976}, which allowed further insight into the method.

ADMM methods are simple to implement and computationally cheap, even for large problems.
However, they tend to converge slowly to a high accuracy solution and the detection of infeasibility is more involved compared to interior-point methods.
They are therefore most often used in applications where a modestly accurate solution is sufficient~\cite{Parikh_2014}.
Most of the early advances in first-order methods such as ADMM happened in the 1970s/80s long before the demand for large scale optimisation, which may explain why they stayed less well-known and have only recently resurfaced.

A method of exploiting the sparsity pattern of PSD constraints in an interior-point algorithm was developed in~\cite{Fukuda_2001}.
This work showed that if the coefficient matrices of an SDP exhibit an aggregate sparsity pattern represented by a chordal graph, then both primal and dual problems can be decomposed into a problem with many smaller PSD constraints on the nonzero blocks of the original matrix variable.
These blocks are associated with subsets, called \emph{cliques}, of graph vertices.
Moreover, it can be advantageous to merge some of these blocks, for example if they overlap significantly.
The optimal way to merge the blocks, or equivalently the corresponding graph cliques, after the initial decomposition depends on the solver algorithm and is still an open question.
Sparsity in SDPs has also been studied for problems with underlying graph structure, e.g.\ optimal power-flow problems in~\cite{Molzahn_2013} and graph optimisation problems in~\cite{Alizadeh_1995}.

\paragraph{Related Work}
Widely used solvers for conic problems, especially SDPs, include \pkg{SeDuMi}~\cite{Sturm_1999}, \pkg{SDPT3}~\cite{Toh_1999} (both open source, MATLAB), and \pkg{MOSEK}~\cite{ApS_2017} (commercial, C) among others.
All of these solvers implement primal-dual interior-point methods.

Both Fukuda~\cite{Fukuda_2001} and Sun~\cite{Sun_2014} developed interior-point solvers that exploit chordal sparsity patterns in PSD constraints.
Some heuristic methods to merge cliques have been proposed for interior-point methods in~\cite{Molzahn_2013} and been implemented in the \pkg{SparseCoLO} package~\cite{Fujisawa_2009} and the \pkg{CHOMPACK} package~\cite{Andersen_2015a}.

Several solvers based on the ADMM method have been released recently.
The solver \pkg{OSQP} \cite{Stellato_2018} is implemented in C and detects infeasibility based on the differences of the iterates~\cite{Banjac_2017}, but only solves LPs and QPs.
The C-based \pkg{SCS}~\cite{ODonoghue_2016} implements an operator splitting method that solves the primal-dual pair of conic programs in order to provide infeasibility certificates.
The underlying homogeneous self-dual embedding method has been extended by~\cite{Zheng_2019} to exploit sparsity and implemented in the MATLAB solver \pkg{CDCS}.
The conic solvers \pkg{SCS} and \pkg{CDCS} are unable to handle quadratic cost functions directly.
Instead they are forced to reformulate problem with quadratic objective functions by adding a second-order cone constraint, which increases the problem size.
Moreover, they rely on primal-dual formulations to detect infeasibility.
\paragraph{Outline}
In Section~\ref{sec:conic_problems} we define the general conic problem format, its dual problem, as well as optimality and infeasibility conditions.
Section~\ref{sec:admm_algorithm} describes the ADMM algorithm that is used by \pkg{COSMO}.
Section~\ref{sec:chordal_decomp} explains how to decompose SDPs in a preprocessing step provided the problem data has an aggregated sparsity pattern.
In Section~\ref{sec:clique_merging} we describe a new clique merging strategy and compare it to existing approaches. Implementation details and code related design choices are discussed in Section~\ref{sec:implementation}.
Section~\ref{sec:results} shows benchmark results of \pkg{COSMO} vs.\ other state-of-the art solvers on a number of test problems. Section~\ref{sec:conclusion} concludes the paper.
\paragraph{Contributions}
With the solver package described in this paper we make the following contributions:
\begin{enumerate}
  \item We implement a first-order method for large conic problems that is able to detect infeasibility without the need of a homogeneous self-dual embedding.
  \item \pkg{COSMO} directly supports quadratic objective functions, thus reducing overheads for applications with both quadratic objective function and PSD constraints.
  This also avoids a major disadvantage of conic solvers compared to native QP solvers, i.e\ no additional matrix factorisation for the conversion is needed and favourable sparsity in the objective can be maintained.
  \item Large structured positive semidefinite programs are analysed and, if possible, chordally decomposed.
  This typically allows one to solve very large sparse and structured problems orders of magnitudes faster than competing solvers.
  For complex sparsity patterns, further performance improvements are achieved by recombining some of the sub-blocks of the initial decomposition in an optimal way.
  For this purpose, we propose a new clique graph based merging strategy and compare it to existing heuristic approaches.
  \item The open-source solver is written in a modular way in the fast and flexible programming language Julia.
  The design allows users to extend the solver by specifying a specific linear system solver and by defining their own convex cones or custom projection methods.
  \end{enumerate}

\paragraph{Notation}
The following notation and definitions will be used throughout this paper.
Denote the space of real numbers $\mb{R}$, the n-dimensional real space $\mb{R}^n$, the n-dimensional zero cone $\{0\}^n$, the nonnegative orthant $\Re_+^n$, the space of symmetric matrices $\mb{S}^n$, and the set of positive semidefinite matrices $\mb{S}^n_+$.

In some of the following sections matrix data is considered in vectorized form.
Denote the vectorization of a matrix $X$ by stacking it columns as $x \eqdef \text{vec}(X)$ and the inverse operation as $\text{vec}^{-1}(X)=\text{mat}(x)$.
For symmetric matrices it is often computationally beneficial to work only with the upper-triangular elements of the matrix.
Denote the transformation of a symmetric matrix $V \in \Sym$ with i,j-th element $V_{ij}$ into a vector of upper-triangular elements as
\begin{equation}
  \mathrm{svec}(V) \eqdef \left[V_{11}, \sqrt{2}V_{12}, V_{22}, \sqrt{2}V_{13}, \sqrt{2}V_{23}, V_{33}, \sqrt{2}V_{14}, \ldots , V_{nn}\right]^\top.
\end{equation}
Here the scaling factor of $\sqrt{2}$ preserves the matrix inner product, i.e.\ $\mathrm{tr}(AB) = \mathrm{svec}(A)^\top\mathrm{svec}(B)$ for symmetric matrices $A,B \in \Sym$.
The inverse operation is denoted by $\mathrm{svec}\inv(s)\reqdef\mathrm{smat}(s)$.
Denote the Kronecker product of two matrices $A$ and $B$ as $ A\otimes B$.
The Frobenius norm of a matrix $A$ is given by $||A||_F=\bigl(\sum_{ij}|a_{ij}|^2)\bigr)^{1/2}$.

Sometimes we consider positive semidefinite constraints in vector form, so we define the space of vectorized symmetric positive semidefinite matrices as
\begin{equation*}
  \mathcal{S}_+^n \eqdef \left\{s \in \mb{R}^{\frac{n(n+1)}{2}} : \text{smat}(s) \in \mb{S}^n_+  \right\}.
\end{equation*}
For a convex cone $\mc{K}$ denote the \emph{polar cone} by
\begin{equation*}
\mc{K}^\circ \eqdef \left\{ y \in \Re^n \mid  \underset{x \in \mc{K}}{\sup} \langle x,y \rangle \leq 0 \right\},
\end{equation*}
the \emph{normal cone} of $\mc{K}$ by
\begin{equation*}
N_\mc{K}(x)  \eqdef \left\{ y \in \Re^n \mid  \underset{\bar{x} \in \mc{K}}{\sup} \langle \bar{x}-x,y \rangle \leq 0 \right\},
\end{equation*}
and, following~\cite{Rockafellar_1970}, the \emph{recession cone} of $\mc{K}$ by
\begin{equation*}
  \mc{K}^{\infty} \eqdef \left\{ y \in \Re^n \mid x + ay \in \mc{K},\, \forall x \in \mc{K}, \forall a \geq 0 \right\}.
\end{equation*}
The \emph{proximal operator} of a convex, closed and proper function $f\colon \mb{R}^n \rightarrow \mb{R}$ is given by
\begin{equation*}
  \prox_f(x) \eqdef \underset{y}{\argmin} \left\{ f(y) + \textstyle{\frac{1}{2}}\norm{y-x}^2_2 \right\}.
\end{equation*}
We denote the \emph{indicator function} of a nonempty, closed convex set $\mc{C} \subseteq \mb{R}^n$ by
\begin{equation*}
  I_{\mc{C}}(x) \eqdef \begin{cases}
    0 & x\in \mc{C}\\
    +\infty & \text{otherwise},
  \end{cases}
\end{equation*}
and the \emph{projection} of $x \in \mb{R}^n$ onto $\mc{C}$ by:
\begin{equation*}
  \project{\mc{C}}(x)\eqdef \underset{y \in \mc{C}}{\argmin} \norm{x-y}_2^2.
\end{equation*}
We further denote the \emph{support function} of $\mc{C}$ by:
\begin{equation*}
  \sigma_{\mc{C}}(x) \eqdef \underset{y \in \mc{C}}{\sup}\langle x,y\rangle.
\end{equation*}


\section{Conic Problems}
\label{sec:conic_problems}
We will address convex optimisation problems with a quadratic objective function and a number of conic constraints in the form:
\begin{equation}
  \begin{array}{ll}
    \mbox{minimize}   & \textstyle{\frac{1}{2}}x^\top Px + q^\top x\\
    \mbox{subject to} & Ax + s  = b \\
                      & s \in \mathcal{K},
  \end{array}
  \label{eq:primal}
\end{equation}
where $x \in \Re^n$ is the primal \emph{decision variable} and $s\in \Re^m$ is the primal \emph{slack variable}.
The objective function is defined by positive semidefinite matrix $P \in \mathbb{S}^n_+$ and vector $q \in \mathbb{R}^n$.
The constraints are defined by matrix $A \in \mathbb{R}^{m\times n}$, vector $b \in \mathbb{R}^m$ and a non-empty, closed, convex cone $\mathcal{K}$ which itself can be a Cartesian product of cones in the form
\begin{equation}
\label{eq:cone_prod}
  \mathcal{K} = \mathcal{K}_1^{m_1} \times \mathcal{K}_2^{m_2} \times \cdots \times \mathcal{K}_N^{m_N},
\end{equation}
with cone dimensions $\sum_{i=1}^N m_i = m$.
Note that any LP, QP, SOCP, or SDP can be written in the form~\eqref{eq:primal} using an appropriate choice of cones, as well as problems involving the power or exponential cones and their duals.

The dual problem associated with~\eqref{eq:primal} is given by:
\begin{alignat}{2}
\label{eq:dual}
&\text{maximize}   & \quad   & -\textstyle{\frac{1}{2}}x^\top Px + b^\top  y - \sigma_\mathcal{K}(y) \notag \\
&\text{subject to} & & Px - A^\top  y = -q\\
&            & & y \in (\mathcal{K}^\infty)^\circ \notag,
\end{alignat}
with dual variable $y \in \mathbb{R}^m$.

The conditions for optimality (assuming linear independence constraint qualification) follow from the Karush-Kuhn-Tucker (KKT) conditions:
\begin{subequations}
\label{eq:optimality}
  \begin{align}
Ax +s &= b, \label{eq:opt1}\\
Px + q - A^\top  y &= 0, \label{eq:opt2}\\
s \in \mathcal{K}, & \quad y \in (\mathcal{K}^\infty)^\circ. \label{eq:opt3}
\end{align}
\end{subequations}
Assuming strong duality, if there exists a $x^* \in \mathbb{R}^n$, $s^* \in \mathbb{R}^m$, and $y^* \in \mathbb{R}^m$ that fulfil~\eqref{eq:opt1}--\eqref{eq:opt3} then the pair $(x^*,s^*)$ is called the primal solution and $y^*$ is called the dual solution of problem~\eqref{eq:primal}.

\subsection{Infeasibility certificates}\label{ssec:infeasibiltyCerts}
Primal and dual infeasibility conditions were developed for ADMM in~\cite{Banjac_2017}.
These conditions are directly applicable to problems of the form~\eqref{eq:primal}.
To simplify the notation of the conditions, define the cone $\bar{\mathcal{K}} \eqdef \mathcal{-K} + \{b\}$.
Then, the following sets provide certificates for primal and dual infeasibility:
\begin{align}
\mathcal{P} &= \left\{x \in \mathbb{R}^n \mid  Px = 0, \, Ax \in \bar{\mathcal{K}}^{\infty}, \, \langle  q,x \rangle < 0  \right\}. \label{eq:d_inf}\\
\mathcal{D} &= \left\{y \in \mathbb{R}^m \mid  A^\top  y  = 0,  \, \sigma_{\bar{\mathcal{K}}}(y) < 0 \right\}, \label{eq:p_inf}
\end{align}
The existence of some $y \in \mathcal{D}$ certifies that problem~\eqref{eq:primal} is primal infeasible, while the existence of some $x \in \mathcal{P}$ certifies dual infeasibility.


\section{ADMM Algorithm}
\label{sec:admm_algorithm}
We use the same splitting as in~\cite{Stellato_2018} to transform problem~\eqref{eq:primal} into standard ADMM format.
The problem is rewritten by introducing the dummy variables $\tx = x $ and $\ts = s$:
\begin{alignat}{2}
\label{eq:splitting}
&\text{minimize}   &  ~   & \textstyle{\frac{1}{2}}\tx^\top P \tx + q^\top \tx + I_{Ax+s=b}(\tx,\ts) + I_{\mathcal{K}}(s)\\
&\text{subject to} & & (\tx,\ts) = (x,s),\notag
\end{alignat}
where the indicator functions of the sets $\{ (x,s) \in \mathbb{R}^n \times \mathbb{R}^m \mid Ax + s = b\}$ and $\mathcal{K}$ were used to move the constraints of~\eqref{eq:primal} into the objective function.
The augmented Lagrangian of~\eqref{eq:splitting} is given by
\begin{equation}
\begin{split}
  L(x,s,\tx,&\ts,\lambda,y) = \textstyle{\frac{1}{2}}\tx^\top P\tx + q^\top \tx + \mathcal{I}_{Ax+s=b}(\tx,\ts) + \mathcal{I}_{\mathcal{K}}(s) \\
&+ \frac{\sigma}{2} \norm{\tx - x + \textstyle{\frac{1}{\sigma}} \lambda}_2^2 + \frac{\rho}{2} \norm{\ts - s + \textstyle{\frac{1}{\rho}} y}_2^2,
\end{split}
\end{equation}
with step size parameters $\rho > 0$ and $\sigma >0$ and dual variables $\lambda \in \mb{R}^n$ and $y \in \mb{R}^m$. The corresponding ADMM iteration is given by:
\begin{subequations}\label{eqn:admmFullAlg}
\begin{align}
    ( \tx^{k+1},&\ts^{k+1})  =\argmin\limits_{\tx,\ts}  L\left( \tx,\ts,x^k,s^k,\lambda^k,y^k \right),\label{eq:min_up}\\
    x^{k+1} &= \alpha \tx^{k+1} + (1-\alpha)x^k + \frac{1}{\sigma} \lambda^k,\label{eq:xup}\\
    s^{k+1} &= \argmin\limits_{s}\frac{\rho}{2} \norm{\alpha \ts^{k+1} + (1-\alpha)s^k - s+\textstyle{\frac{1}{\rho}}y^k}_2^2 + I_{\mathcal{K}}(s), \label{eq:sup}\\
    \lambda^{k+1} &= \lambda^k + \sigma \left(\alpha \tx^{k+1} + (1-\alpha)x^{k} - x^{k+1}\right),\label{eq:lup}\\
    y^{k+1} &= y^k + \rho \left(\alpha \ts^{k+1} +(1-\alpha)s^k -s^{k+1} \right),\label{eq:up2}
\end{align}
\end{subequations}
where we relaxed the $z$-update and the dual variable update with relaxation parameter $\alpha \in (0,2)$ according to~\cite{Eckstein_1992}.
Notice from~\eqref{eq:xup} and~\eqref{eq:lup} that the dual variable corresponding to the constraint $x=\tx$ satisfies $\lambda^k=0 \text{ for all } k$.

\subsection{Solution of the equality constrained QP}
The minimization problem in~\eqref{eq:min_up} has the form of an equality-constrained quadratic program:
\begin{alignat}{2}
\label{eq:ec_qp}
&\text{minimize}  & \,    & \textstyle{\frac{1}{2}}\tx^\top P\tx + q^\top \tx + \frac{\sigma}{2} \norm{\tx - x^k }_2^2 + \frac{\rho}{2} \norm{\ts - s^k + \textstyle{\frac{1}{\rho}} y^k}_2^2\\
&\text{subject to } & & A\tx+\ts=b.\notag
\end{alignat}
The solution of~\eqref{eq:ec_qp} can be obtained by solving a single linear system.
The corresponding Lagrangian is given by:
\begin{equation}
\mathcal{L}(\tx,\ts,\nu) =\textstyle{\frac{1}{2}}\tx^\top P\tx + q^\top \tx + \frac{\sigma}{2} \norm{\tx - x^k}_2^2 + \frac{\rho}{2} \norm{\ts - s^k + \frac{1}{\rho} y^k}_2^2+ \nu^\top\left(A\tx+\ts-b\right),
\end{equation}
where the Lagrangian multiplier $\nu \in \mathbb{R}^m$ accounts for the equality-constraint $Ax+s = b$. Thus, the KKT optimality conditions for this equality constrained QP are given by:
    \begin{align}
\frac{\partial \mathcal{L}}{\partial \tx} = P\tx^{k+1} + q +\sigma \left(\tx^{k+1}-x^k\right) + A^\top\nu^{k+1} &= 0,\\
\frac{\partial \mathcal{L}}{\partial \ts} = \rho \left(\ts^{k+1}-s^k+\frac{1}{\rho} y^k\right) + \nu^{k+1} &= 0,\\
A\tx^{k+1} +\ts^{k+1} - b &= 0.
\end{align}
Elimination of $\ts^{k+1}$ from these equations leads to the linear system:
\begin{align}\label{eq:sys}
\begin{bmatrix}
P + \sigma I & A^\top \\A &- \frac{1}{\rho}I
    \end{bmatrix}\begin{bmatrix}\tx^{k+1} \\ \nu^{k+1}\end{bmatrix}&= \begin{bmatrix}-q+\sigma x^k \\b-s^k+\frac{1}{\rho}y^k\end{bmatrix}
\end{align}
with
\begin{equation}
\ts^{k+1} = s^k - \frac{1}{\rho}\left(\nu^{k+1} + y^{k}\right).\label{eq:seq}
\end{equation}
Note that the introduction of the dummy variable $\tx$ led to the term $\sigma I$ in the upper-left corner of the coefficient matrix in~\eqref{eq:sys}.
Consequently, the coefficient matrix in \eqref{eq:sys} is always quasi-definite~\cite{Vanderbei_1995}, i.e.\ it always has a positive definite upper-left block and a negative definite lower-right block, and is therefore full rank even when $P = 0$ or $A$ is rank deficient.
Following~\cite{Vanderbei_1995} the left hand side of \eqref{eq:sys} always has a well-defined $LDL^\top$ factorization with a diagonal $D$.
\subsection{Projection step}
As shown in~\cite{Parikh_2014} the minimization problem in~\eqref{eq:sup} can be interpreted as the $\rho$-weighted proximal operator of the indicator function $I_\mathcal{K}$.
It is therefore equivalent to the Euclidean projection $\Pi_\mc{K}$ onto the cone $\mathcal{K}$, i.e.\
  \begin{equation}
  \label{eq:proj_step}
  s^{k+1} = \Pi_{\mathcal{K}}\left(\alpha\ts^{k+1}+(1-\alpha)s^k+\frac{1}{\rho}y^k\right).
  \end{equation}
If $\mc{K}$ is a Cartesian product of cones as in~\eqref{eq:cone_prod} this projection is equivalent to the projection of the relevant components of the argument of $\Pi_\mc{K}(\cdot)$ onto each cone $\mc{K}_i$.
A problem with $N$ SDP constraints therefore requires $N$ projections, but, since each of these operates on an independent segment of the input vector, they can be performed in parallel.

\subsection{Algorithm steps}
\label{ssec:algorithm}
The calculations performed at each iteration are summarized in Algorithm~\ref{alg:admm}.
\begin{figure}[tb]
\begin{algorithm}[H]
\SetKwInOut{Input}{Input}
\SetKwInOut{Output}{Output}
\SetKwRepeat{Do}{Do}{while}

\Input{ initial values $x^0$, $s^0$, $y^0$, problem data $P$, $q$, $A$, $b$, and parameters $\sigma > 0$, $\rho > 0$, $\alpha \in (0,2)$}
 \Do{termination criteria not satisfied}{
  $(\tx^{k+1},\nu^{k+1}) \, \leftarrow$ solve linear system \eqref{eq:sys}\;
  $\ts^{k+1}  \, \leftarrow \, s^k - \frac{1}{\rho}(\nu^{k+1} + y^{k})$\;
  $ x^{k+1} \, \leftarrow \, \alpha \tx^{k+1} + (1-\alpha)x^k$\;
  $ s^{k+1} \, \leftarrow \, \Pi_{\mathcal{K}}\left(\alpha\ts^{k+1}+(1-\alpha)s^k+\frac{1}{\rho}y^k\right)$\;
  $y^{k+1} \, \leftarrow \,y^k + \rho (\alpha \ts^{k+1} +(1-\alpha)s^k -s^{k+1} )$\;
 }
 \caption{ADMM iteration}
   \label{alg:admm}
\end{algorithm}
\end{figure}
Observe that the coefficient matrix of the linear system in~\eqref{eq:sys} is constant, so that one can precompute and cache the $LDL^\top$ factorization and efficiently evaluate line 2 with changing right hand sides.
A refactorisation is only necessary if the step size parameters $\rho$ and $\sigma$ are updated.

Lines 3, 4, and 6 are computationally inexpensive since they involve only vector addition and scalar-vector multiplication.
The projection in line 5 is crucial to the performance of the algorithm depending on the particular cones employed in the model: projections onto the zero-cone or the nonnegative orthant are inexpensive, while a projection onto the positive-semidefinite cone of dimension $N$ involves an eigenvalue decomposition.
Since direct methods for eigen-decompositions have a complexity of approximately $\mathcal{O}(N^3)$, this turns line 5 into the  most computationally expensive operation of the algorithm for large SDPs, and improving the efficiency of this step will be the objective of much of Sections~\ref{sec:chordal_decomp} and~\ref{sec:clique_merging}.

\subsection{Algorithm convergence}\label{ssec:algorithmConvergence}
For feasible problems, Algorithm~\ref{alg:admm} produces a sequence of iterates $(x^k,s^k,y^k)$ that converges to a limit satisfying the optimality conditions in~\eqref{eq:optimality} as $k\rightarrow \infty$.
The authors in~\cite{Stellato_2018} show convergence of Algorithm~\ref{alg:admm} by applying the Douglas-Rachford splitting to a problem reformulation.

In the Douglas-Rachford formulation, lines 5 and 6 become projections onto $\mc{K}$ and $\mc{K}^\circ$ respectively, so that the conic constraints in~\eqref{eq:opt3} always hold.
Furthermore, convergence of the residual iterates in~\eqref{eq:opt1}--\eqref{eq:opt2} can be concluded from the convergence of the splitting variables
\begin{equation}
  x^k - \tx^k \rightarrow 0,\quad s^k-\ts^k\rightarrow 0,
\end{equation}
which generally holds for Douglas-Rachford splitting~\cite{Bauschke_2011a}.

For infeasible problems,~\cite{Banjac_2017} showed that Algorithm~\ref{alg:admm} leads to convergence of the successive differences between iterates
\begin{equation}
  \delta x^k = x^k - x^{k-1},\, \delta s^k = s^k - s^{k-1},\, \delta y^k = y^k - y^{k-1}.
\end{equation}
For primal infeasible problems $\delta y = \lim_{k\rightarrow \infty}\delta y^k$ will satisfy condition~\eqref{eq:p_inf}, whereas for dual infeasible problems $\delta x = \lim_{k\rightarrow \infty}\delta x^k$ is a certificate of~\eqref{eq:d_inf}.

\subsection{Scaling the problem data}
The rate of convergence of ADMM and other first-order methods depends in practice on the scaling of the problem data; see~\cite{Giselsson_2014}.
Particularly for badly conditioned problems, this suggests a preprocessing step where the problem data is scaled in order to improve convergence.
For certain problem classes an optimal scaling has been found, see~\cite{Giselsson_2014,Giselsson_2015,Greenbaum_1997}.
However, the computation of the optimal scaling is often more complicated than solving the original problem.
Consequently, most algorithms rely on heuristic methods such as matrix equilibration.

We scale the equality constraints by diagonal positive definite matrices $D$ and $U$. The scaled form of~\eqref{eq:primal} is given by:
\begin{alignat}{2}
\label{eq:scaled}
&\text{minimize}   & \quad   & \textstyle{\frac{1}{2}} \hx^\top \hat{P} \hx + \hat{q}^\top \hx\\
&\text{subject to} & & \hat{A} \hx + \hat{s}  = \hat{b}, \notag \\
&            & & \hs \in U\mathcal{K}\notag,
\end{alignat}
with scaled problem data
\begin{equation}
\label{eq:scaling1}
  \hat{P}=DPD, \quad \hat{q}=Dq,  \quad\hat{A}=UAD, \quad \hat{b}=Ub,
\end{equation}
and the scaled convex cone $U\mc{K} \eqdef \{Uv \in \mb{R}^m \mid v \in \mc{K} \}$.
After solving~\eqref{eq:scaled} the original solution is obtained by reversing the scaling:
\begin{equation}
\label{eq:scaling2}
   x = D\hx, \quad s = U^{-1}\hs, \quad y = U\hat{y}.
 \end{equation}
One heuristic strategy that has been shown to work well in practice is to choose the scaling matrices $D$ and $U$ to equilibrate, i.e.\ reduce the condition number of, the problem data.
The Ruiz equilibration technique described in~\cite{Ruiz_2001} iteratively scales the rows and columns of a matrix to have an infinity-norm of $1$ and converges linearly.
We apply the modified Ruiz algorithm shown in Algorithm~\ref{alg:ruiz} to reduce the condition number of the symmetric matrix $R$
\begin{equation}
 R =\begin{bmatrix}P & A^\top\\A & 0\end{bmatrix}
\end{equation}
which represents the problem data.
\begin{figure}[h]
\begin{algorithm}[H]
\SetKw{Set}{set}
\SetKw{Return}{return}
\SetKwInOut{Output}{Output}
\SetKwRepeat{Do}{Do}{while}

\Set{$D=I_\text{n}$, $U=I_\text{m}$, $c=\mathbf{1}_\text{m+n}$}\;
 \Do{$\norm{\mathbf{1}-c}_\infty > \mathrm{tol}$}{
  \For{$i=1,\ldots,n+m$}{
    \If{$\norm{R_{c,i}}_\infty > \tau$}{
    $c_i \leftarrow \norm{R_{c,i}}_\infty^{-\frac{1}{2}}$;
    }
    }
    $\hat{D} = \text{diag}(c_{1:n})$, \, $\hat{U} = \text{diag}(c_{n+1:n+m})$\;
    $D = \hat{D} \cdot D$, \, $U = \hat{U} \cdot U$\;
    $P = \hat{D}P\hat{D}$, \, $A = \hat{U}A\hat{D}$\;
    assemble $R$\;
 }
 \Return{$D$, $U$}\;
 \caption{Modified Ruiz equilibration}
   \label{alg:ruiz}
\end{algorithm}
\end{figure}
Since $R$ is symmetric it suffices to consider the columns $R_{c,i}$ of $R$.
At each iteration the scaling routine calculates the norm of each column.
For the columns with norms higher than the tolerance $\tau$ scaling vector $c$ is updated with the inverse square root of the norm\footnote{For the presented results a value of $\tau=10^{-6}$ was chosen. }.
If the norm is below the tolerance, the corresponding column will be scaled by 1.

Since the matrix $U$ scales the (possibly composite) cone constraint, the scaling must ensure that if $s\in \mc{K}$ then \mbox{$U^{-1}s \in \mc{K}$}.
Let $\mc{K}$ be a Cartesian product of $N$ cones as in~\eqref{eq:cone_prod} and partition $U$ into blocks
\begin{equation}
  U=\text{diag}\left(U_1,\ldots,U_N\right),
\end{equation}
with block $U_i \in \mb{R}^{m_i \times m_i}$ which scales the constraint corresponding to $\mc{K}_i$.
For each cone $\mc{K}_i \in \mathbb{R}^{m_i}$ that requires a scalar or symmetric scaling, \eg a second-order cone or positive semidefinite cone, the corresponding block $U_i$ is replaced with
\begin{equation}
  U_i^* \eqdef u_i I_{m_i},\quad \text{for } i=1,\ldots,N
\end{equation}
where the mean value of the diagonal entries of the original block in $U$, $u_i= \text{tr}(U_i)/m_i$, was chosen as a heuristic scaling factor.

\subsection{Termination criteria}
\label{sec:termination}
The termination criteria discussed in this section are based on the unscaled problem data and iterates.
Thus, before checking for termination the solver first reverses the scaling according to equations~\eqref{eq:scaling1}--\eqref{eq:scaling2}.
To measure the progress of the algorithm, we define the primal and dual residuals of the problem as:
\begin{subequations}
\label{eq:residuals}
\begin{align}
r_p &\eqdef Ax + s -b,\label{eq:rp}\\
r_d &\eqdef Px + q - A^\top y.\label{eq:rd}
\end{align}
\end{subequations}
According to Section 3.3 of~\cite{Boyd_2011} a valid termination criterion is that the size of the norms of the residual iterates in~\eqref{eq:residuals} are small.
Our algorithm terminates if the residual norms are below the sum of an absolute and a relative tolerance term:
\begin{subequations}
\begin{align}
  \norm{r_p^k}_\infty &\leq \epsilon_{\mathrm{abs}} + \epsilon_{\mathrm{rel}} \, \text{max} \left\{ \norm{Ax^k}_\infty,\norm{s^k}_\infty, \norm{b}_\infty \right\},\\
   \norm{r_d^k}_\infty &\leq \epsilon_{\mathrm{abs}} + \epsilon_{\mathrm{rel}} \, \text{max} \left\{\norm{Px^k}_\infty,\norm{q}_\infty, \norm{A^\top y^k}_\infty \right\},
\end{align}
\end{subequations}
where $\epsilon_{\mathrm{abs}}$ and $\epsilon_{\mathrm{rel}}$ are user defined tolerances.

Following~\cite{Banjac_2017}, the algorithm determines if the one-step differences $\delta x^k$ and $\delta y^k$ of the primal and dual variable fulfil the normalized infeasibility conditions~\eqref{eq:p_inf}--\eqref{eq:d_inf} up to certain tolerances $\epsilon_{\text{p,inf}}$ and $\epsilon_{\mathrm{d,inf}}$.
The solver returns a primal infeasibility certificate if
\begin{subequations}
\begin{align}
  \norm{A^\top \delta y^k}_\infty / \norm{\delta y^k}_\infty &\leq \epsilon_{\mathrm{p,inf}}, \\
  \sigma_{\bar{\mathcal{K}}}\left(\delta y^k\right) &\leq \epsilon_{\mathrm{p,inf}},
\end{align}
\end{subequations}
holds and a dual infeasibility certificate if
\begin{subequations}
  \begin{align}
  \norm{P \delta x^k}_\infty / \norm{\delta x^k}_\infty &\leq \epsilon_{\mathrm{d,inf}},  \\
  q^\top \delta x^k / \norm{\delta x^k}_\infty &\leq \epsilon_{\mathrm{d,inf}},  \\
   A \delta x^k + v &\in \bar{\mathcal{K}}^\infty, \label{eq:dinf_cond}\\
  \text{with } \norm{v}_\infty &\leq \epsilon_{\mathrm{d,inf}} \norm{\delta x^k}_\infty, \notag
\end{align}
\end{subequations}
holds.


\section{Chordal Decomposition}
\label{sec:chordal_decomp}
As noted in~Section~\ref{ssec:algorithm}, for large SDPs the eigen-decomposition in the projection step~\eqref{eq:proj_step} is the principal performance bottleneck for the algorithm.
However, since large-scale SDPs often exhibit a certain structure or sparsity pattern, a sensible strategy is to exploit any such structure to alleviate this bottleneck.
If the aggregated sparsity pattern is \emph{chordal}, Agler's~\cite{Agler_1988} and Grone's~\cite{Grone_1984} theorems can be used to decompose a large PSD constraint into a collection of smaller PSD constraints and additional coupling constraints.
The projection step applied to the set of smaller PSD constraints is usually significantly faster than when applied to the original constraint.
Since the projections are independent of each other, further performance improvement can be achieved by carrying them out in parallel.
Our approach is to apply chordal decomposition to a standard form SDP in the form \eqref{eq:primal} to produce a decomposed problem, and then transform the resulting decomposed problem back to a  problem in the form \eqref{eq:primal} but with more variables and a collection of smaller PSD cone constraints.
This process allows us to apply chordal decomposition as a preprocessing step before the problem is handed to the solver.
As discussed in~Section~\ref{sub:agler_s_theorem}, a chordal sparsity pattern can be imposed on any PSD constraint.

\subsection{Graph preliminaries} 
\label{sub:graph_preliminaries}
In the following we define graph-related concepts that are useful to describe the sparsity structure of a problem.
A good overview of this topic is given in the survey paper~\cite{Vandenberghe_2015}.
Consider the \emph{undirected graph} $G(V,E)$ with vertex set $V=\{1,\ldots,n\}$ and edge set $E \subseteq V \times V$.
If all vertices are pairwise adjacent, i.e.\ $E = \left\{ \{v,u\} \mid v,u \in V, \, v \neq u \right\}$, the graph is called \emph{complete}.
We follow the convention of~\cite{Vandenberghe_2015} by defining a \emph{clique} as a subset of vertices $\clique{} \subseteq V$ that induces a \emph{maximal} complete subgraph of $G$.
The number of vertices in a \emph{clique} is given by the cardinality $|\clique{}|$.
A \emph{cycle} is a path of edges (i.e.\ a sequence of distinct edges) joining a sequence of vertices in which only the first and last vertices are repeated.

The following decomposition theory relies on a subset of graphs that exhibit the important property of \emph{chordality}.
A graph $G$ is \emph{chordal} (or \emph{triangulated}) if every cycle of length greater than three has a \emph{chord}, which is an edge between nonconsecutive vertices of the cycle.
A non-chordal graph can always be made chordal by adding extra edges.

\begin{figure*}[htb]
\centering
\begin{tabular}[b]{@{}cccc@{}}
\begin{minipage}{0.3\textwidth}
\centering
\includegraphics[width = 0.85\textwidth]{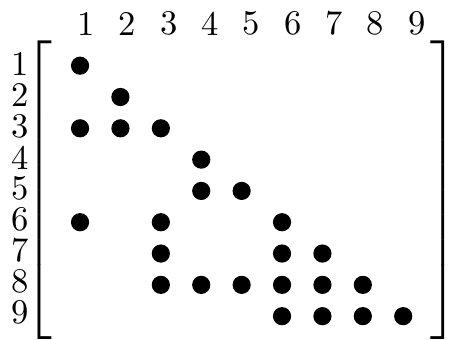}
\end{minipage}
&
\begin{minipage}{0.3\textwidth}
\centering
\includegraphics[width = 0.85\textwidth]{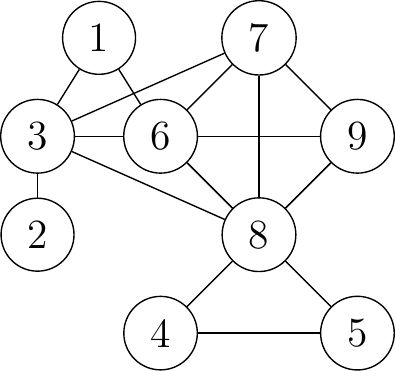}

 \end{minipage}
&
\begin{minipage}{0.3\textwidth}
\centering
\includegraphics[width = 0.85\textwidth]{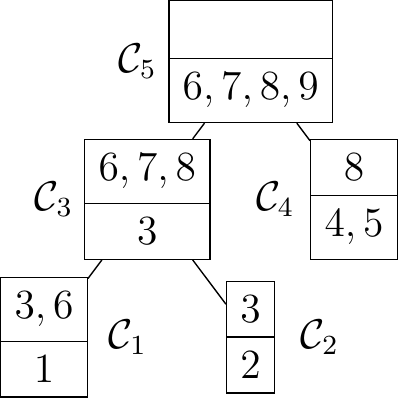}

\end{minipage}
\\
(a) & (b) & (c)
\end{tabular}

  \caption{(a)~Aggregate sparsity pattern, (b)~sparsity graph $G(V,E)$, and (c)~clique tree $\mc{T}(\mc{B}, \mc{E})$.}
  \label{fig:example_structures}
\end{figure*}

An undirected graph with $n$ vertices can be used to represent the sparsity pattern of a symmetric matrix $S \in \Sym[n]$.
Every nonzero entry $S_{ij} \neq 0$ in the lower triangular part of the matrix introduces an edge $(i,j)\in E$.
An example of a sparsity pattern and the associated graph is shown in Figure~\ref{fig:example_structures}(a--b).

For a given sparsity pattern $G(V, E)$, we define the following symmetric sparse matrix cones:
\begin{align}
\label{eq:sparsecone}
\Sym\left(E, 0\right) &\eqdef \left\{ S \in \Sym \mid S_{ij} = S_{ji} = 0 \; \text{if } i \neq j \text{ and } (i,j) \notin E\right\},\\
\Psd\left(E, 0\right) &\eqdef \left\{ S \in \Sym(E, 0) \mid S \succeq 0\right\}.
\end{align}
Given the definition in~\eqref{eq:sparsecone} and a matrix $S \in \Sym\left(E, 0\right)$, note that the diagonal entries $S_{ii}$ and the off-diagonal entries $S_{ij}$ with $(i,j)\in E$ may be zero or nonzero.
Moreover, we define the cone of positive semidefinite completable matrices as
\begin{equation}
  \Psd(E, ?) \eqdef \left\{ Y \mid \exists  \hat{Y} \in \Psd,  Y_{ij} = \hat{Y}_{ij}, \text{if}\, i=j \,\text{or}\, (i,j) \in \! E \! \right\}.
\end{equation}
For a matrix $Y \in  \Psd(E, ?)$ we can find a positive semidefinite completion by choosing appropriate values for all entries $(i, j)  \notin E$.
An algorithm to find this completion is described in~\cite{Vandenberghe_2015}.
An important structure that conveys a lot of information about the nonzero blocks of a matrix, or equivalently the cliques of a chordal graph, is the \emph{clique tree} (or \emph{junction tree}).
For a chordal graph $G$ let $\mc{B} = \{\mc{C}_1, \ldots, \mc{C}_p \}$ be the set of cliques.
The clique tree $\mc{T}(\mc{B},\mc{E})$ is formed by taking the cliques as vertices and by choosing edges from $\mc{E} \subseteq \mc{B} \times \mc{B}$ such that the tree satisfies the \emph{running-intersection property}:

\begin{definition}[Running intersection property]~\\
For each pair of cliques $\clique{i}$, $\clique{j} \in \mc{B}$, the intersection $\clique{i} \cap \clique{j}$ is contained in all the cliques on the path in the clique tree connecting $\clique{i}$ and $\clique{j}$.
\end{definition}
This property is also referred to as the \textit{clique-intersection property} in~\cite{Nakata_2003} and the \textit{induced subtree property} in~\cite{Vandenberghe_2015}.
For a given chordal graph, a clique tree can be computed using the algorithm described in~\cite{Pothen_1990}.
The clique tree for an example sparsity pattern is shown in Figure~\ref{fig:example_structures}(c).

In a \emph{post-ordered} clique tree the descendants of a node are given consecutive numbers, and
a suitable post-ordering can be found via depth-first search.
For a clique $\clique{\ell}$ we refer to the first clique encountered on the path to the root as its \textit{parent clique} $\clique{\mathrm{par}}$.
Conversely $\clique{\ell}$ is called the \textit{child} of  $\clique{\mathrm{par}}$.
If two cliques have the same parent clique we refer to them as \textit{siblings}.

We define the function $\mathrm{par} : 2^V \to 2^V$ and the multivalued function $\mathrm{ch} : 2^V \rightrightarrows 2^V$ such that $\mathrm{par}(\clique{\ell})$ and $\mathrm{ch}(\clique{\ell})$ return the parent clique and set of child cliques of $\clique{\ell}$, respectively, where $2^{V}$ is the power set (set of all subsets) of $V$.

Note that each clique in Figure~\ref{fig:example_structures}(c) has been partitioned into two sets.
The upper row represents the \textit{separators} $\eta_\ell = \clique{\ell} \cap \mathrm{par}(\clique{\ell})$, i.e.\ all clique elements that are also contained in the parent clique.
We call the sets of the remaining vertices shown in the lower rows the \textit{clique residuals} or \textit{supernodes} \mbox{$\nu_\ell = \clique{\ell} \setminus \eta_\ell$}.
Keeping track of which vertices in a clique belong to the supernode and the separator is useful as the information is needed to perform a positive semidefinite completion.
Following the authors in~\cite{Habib_2012}, we say that two cliques $\clique{i}$, $\clique{j}$ form a separating pair $\mc{P}_{ij} = \{ \clique{i},\clique{j} \}$ if their intersection is non-empty and if every path in the underlying graph $G$ from a vertex $\clique{i} \setminus \clique{j}$ to a vertex $\clique{j} \setminus \clique{i}$ contains a vertex of the intersection $\clique{i} \cap \clique{j}$.

\subsection{Agler's and Grone's theorems} 
\label{sub:agler_s_theorem}
To explain how the concepts in the previous section can be used to decompose a positive semidefinite constraint, we first consider an SDP in standard primal form:
\begin{equation}
\OptProblem[eq:primal_sdp]{minimize}{}{\langle C, X \rangle}{ \langle A_k, X \rangle = b_k, \quad k =1,\ldots, m\\ & X \in \Psd,}
\end{equation}
with variable $X$ and coefficient matrices $A_k, C \in \Sym$. The corresponding dual problem is
\begin{equation}
\OptProblem[eq:dual_sdp]{maximize}{}{b^\top y}{\displaystyle \sum_{k = 1}^{m}A_k y_k + S  = C\\ & S \in \Psd,}
\end{equation}
with dual variable $y \in \Re^m$ and slack variable $S$.
Let us assume that the problem matrices in~\eqref{eq:primal_sdp} and~\eqref{eq:dual_sdp} each have their own sparsity pattern
\begin{equation*}
  A_k \in \Sym[n](E_{A_k}, 0) \text{ and } C \in \Sym[n](E_{C}, 0).
\end{equation*}
The \emph{aggregate sparsity} of the problem is given by the graph $G(V, E)$ with edge set
\begin{equation*}
  E = E_{A_1} \cup E_{A_2} \cup \cdots \cup E_{A_m} \cup E_{C}.
\end{equation*}
In general $G(V, E)$ will not be chordal, but a \emph{chordal extension} can be found by adding edges to the graph.
We denote the extended graph as $G(V,\bar{E})$, where $E \subseteq \bar{E} $.
Finding the minimum number of additional edges required to make the graph chordal is an NP-complete problem~\cite{Yannakakis_1981}.
Consider a $0$--$1$ matrix $M$ with edge set $E$.
A commonly used heuristic method to compute the chordal extension is to perform a symbolic Cholesky factorisation of $M$~\cite{Fukuda_2001}, with the edge set of the Cholesky factor then providing the chordal extension $\bar{E}$.
To simplify the notation in the remainder of the article, we henceforward assume that $E$ represents a chordal graph or has been appropriately extended.

Using the aggregate sparsity of the problem we can modify the constraints on the matrix variables in~\eqref{eq:primal_sdp} and~\eqref{eq:dual_sdp} to be in the respective sparse positive semidefinite matrix spaces:
\begin{equation}
\label{eq:sparse_constraints}
  X \in \Psd(E, ?) \; \text{and} \;  S \in \Psd(E, 0).
\end{equation}
We further define the entry-selector matrices \mbox{$T_\ell \in \Re^{ |\clique{\ell}|  \times n }$} for a clique $\clique{\ell}$ as
\begin{equation}
  (T_\ell)_{ij} \eqdef \begin{cases}
    1, & \text{if } \clique{\ell}(i) = j\\
    0,  & \mathrm{otherwise,}
  \end{cases}
\end{equation}
where $\clique{\ell}(i)$ is the $i$th vertex of $\clique{\ell}$.
We can express the constraints in~\eqref{eq:sparse_constraints} in terms of multiple smaller coupled constraints using Grone's and Agler's theorems.

\begin{theorem}[Grone's theorem~\cite{Grone_1984}]\label{thm:grone}
Let $G(V,E)$ be a chordal graph with a set of maximal cliques $\{ \clique{1},\ldots,\clique{p}\}$. Then $X \in \Psd(E, ?)$ if and only if
\begin{equation}
  X_\ell = T_\ell X T_\ell^\top \in \Psd[|\clique{\ell}|],
\end{equation}
for all $\ell=1,\ldots, p$.
\end{theorem}
Applying this theorem to~\eqref{eq:primal_sdp} while restricting $X$ to the positive semidefinite completable matrix cone as in~\eqref{eq:sparse_constraints} yields the decomposed problem:
\begin{equation}
\OptProblemExtra[eq:primal_decomp_sdp]{minimize}{}{\langle C, X \rangle}{ \langle A_k, X \rangle = b_k, & k =1,\ldots, m\\ & X_\ell = T_\ell X T^\top_\ell, & \ell = 1,\ldots,p \\ & X_\ell \in \Psd[|\mc{C}_\ell|], & \ell = 1,\ldots, p,}
\end{equation}

For the dual problem we utilise Agler's theorem, which is the dual to Theorem~\ref{thm:grone}:
\begin{theorem}[Agler's theorem~\cite{Agler_1988}]\label{thm:agler}
Let $G(V,E)$ be a chordal graph with a set of maximal cliques $\{ \clique{1},\ldots,\clique{p}\}$.
Then $S \in \mathbb{S}^n_+(E,0)$ if and only if there exist matrices $S_\ell \in \Psd[|\clique{\ell}|]$ for $\ell=1,\ldots,p$ such that
\begin{equation}
    S = \sum_{\ell=1}^p T_\ell^\top S_\ell T_\ell. \label{eq:agler}
\end{equation}
\end{theorem}

Note that the matrices $T_\ell$ serve to extract the submatrix $S_\ell$ such that $S_\ell=T_\ell S T_\ell^\top $ has rows and columns corresponding to the vertices of the clique $C_\ell$.
With this theorem, we transform the dual form SDP in~\eqref{eq:dual_sdp} with the restriction on $S$ in~\eqref{eq:sparse_constraints} to arrive at:
\begin{equation}
\OptProblem[eq:dual_decomp_sdp]{maximize}{}{ b^\top y}{
\displaystyle \sum_{k = 1}^{m}A_k y_k + \sum_{\ell = 1}^p  T^\top_\ell S_\ell T_\ell = C\\
  & S_\ell \in   \Psd[|\mc{C}_\ell|],  \quad \ell = 1,\ldots, p.}
\end{equation}

\subsection{Returning a decomposed problem into standard SDP form}
After the decomposition results of~Section~\ref{sub:agler_s_theorem} have been applied, the SDP problem~\eqref{eq:dual_decomp_sdp} has to be transformed back into standard form~\eqref{eq:primal}.
For the undecomposed problem in~\eqref{eq:dual_sdp}, this is achieved by first relabeling $y$ as $x$. We then choose $P = 0_{n \times n}$, $q = -b$, define $A = [\svec{A_1}, \ldots, \svec{A_m} ]$, $b = \svec{C}$ and $s = \svec{S}$.

To transform the decomposed dual problem~\eqref{eq:dual_decomp_sdp}, we will make use of the fact that the decision variable $S \in \Sym$ and all the submatrices $S_\ell$ are symmetric and consider instead $s_\ell = \mathrm{svec}(S_\ell)$.
The main challenge is to keep track of the overlapping entries in the individual blocks and ensure they sum to the original entry in $S$. Different possible transformations achieving this are described in~\cite{Kim_2011}.

All the necessary information about the overlapping entries is stored in the clique tree~$\mc{T}(\mc{B}, \mc{E})$ that represents the sparsity pattern of $S$.
We assume that the clique tree is post-ordered with cliques $\clique{1},\ldots,\clique{p}$.
Define a vector to represent slack variables for overlapping entries $\theta \eqdef [\theta_1, \ldots, \theta_{n_o}]^\top$, where $n_o \eqdef \sum_{\ell=1}^p \frac{\abs{\eta_\ell} (\abs{\eta_\ell} + 1)}{2}$ is the total number of overlapping entries in the upper triangle of the sparsity pattern.
The $s$ vector of the decomposed problem is created by stacking the vectorized subblocks $s_\ell$ according to the postordering of the clique tree.
Define $\omega(i,j,\ell)$ as the index of $s$ corresponding to the $(i,j)$th element of block $S_\ell$.
The equality constraint in problem~\eqref{eq:dual_decomp_sdp} can be represented in the equivalent standard form \eqref{eq:primal} as:
 \begin{equation}\label{eq:decomp_transformation}
  \left[
  \begin{array}{c c c | c}
  \svec{\tilde{A}_{1}^{1}} & \cdots & \svec{\tilde{A}_{m}^{1}} & \\
  \vdots & \ddots & \vdots & L\\
  \svec{\tilde{A}_{1}^{p}} & \cdots & \svec{\tilde{A}_{m}^{p}} &
  \end{array}
  \right]
  \begin{bmatrix}
  x \\ \theta
  \end{bmatrix}
  +
  \begin{bmatrix}
    s_1 \\ s_2 \\ \vdots \\ s_p
  \end{bmatrix}
   =
  \begin{bmatrix}
  \svec{\tilde{B}^{1}}  \\ \vdots \\ \svec{\tilde{B}^{p}}
  \end{bmatrix},
 \end{equation}
with
\begin{align}
  \tilde{A}_{k, ij}^\ell \eqdef
  \begin{cases}
  A_{k, ij}^\ell & \text{if } (i, j) \notin \sep{\clique{\ell}}\\
  0         & \text{otherwise},
  \end{cases}
  &&
  \tilde{B}_{ij}^\ell \eqdef
  \begin{cases}
  B_{ij}^\ell & \text{if } (i, j) \notin \sep{\clique{\ell}}\\\
  0         & \text{otherwise}.
  \end{cases}
\end{align}
The matrices $\tilde{A}_{k}^\ell$ and $\tilde{B}^\ell$ take on the values of $A_k$ and $B$ in the locations corresponding to the elements in the submatrix $S_\ell$.
If a matrix entry of clique $\clique{l}$ in $A_{k, ij}^\ell$ and $B_{ij}^\ell$ is overlapped by an entry of the parent clique, i.e.\  both $(i,j)$ are included in $\sep{\clique{l}}$, it is set to zero.
Each column of the linking matrix $L \in \Re^{m_d \times n_o}$ links one overlapping entry in the clique tree.
$L$ is generated by first collecting all the matrix indices $(i, j)_\ell$ of the separators $\sep{\clique{\ell}}$ in the clique tree:
\begin{equation}
  O \eqdef \bigcup_{\ell = 1}^p \left\{ (i, j)_\ell \in \sep{\clique{\ell}} \mid i\leq j\right\}.
\end{equation}
Then $L$ is constructed column-by-column, each representing one overlapping entry.
The column vector $l_c$ is equal to $1$ in the row $r$ corresponding to the $(i,j)_\ell$-th entry of $S_\ell$ and $-1$ in the row corresponding to the same entry of the parent block $S_k$, where $\clique{\ell} = \mathrm{ch}(\clique{k})$.
Thus, each element in $l_c$ is defined by:
\begin{equation}
  l_{rc} \eqdef \begin{cases}
  \hphantom{-}1 & \text{if } r = \omega(i,j, \ell)\\
  -1 & \text{if } r = \omega(i,j, \mathrm{par}(\ell))\\
  \hphantom{-}0 & \text{otherwise},
  \end{cases} \quad \text{for each } (i,j)_\ell \in O.
\end{equation}

As an example consider a problem with $m=1$ and $p=3$ and the simple sparsity pattern and clique tree given by:

 \begin{minipage}{0.8\textwidth}
 \includegraphics{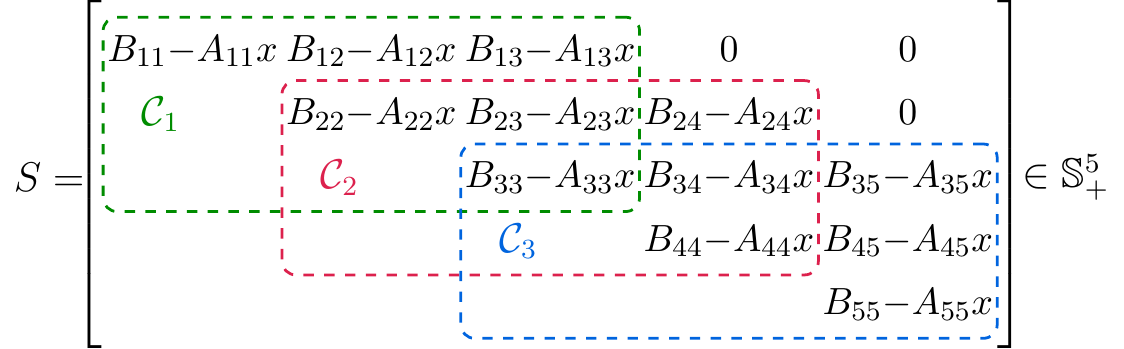}

\end{minipage}%
\begin{minipage}{0.05\textwidth}
    \begin{flushright}
\includegraphics{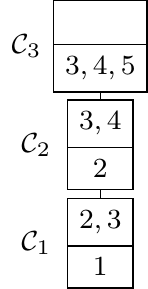}


\end{flushright}
\end{minipage}

Using the transformation in~\eqref{eq:decomp_transformation} this constraint can be represented by three constraints on the submatrices $S_1$, $S_2$ and $S_3$ and six overlap variables $\theta_1,\dots,\theta_6$:
\begin{align*}
 S_3 &= \begin{bmatrix}
  B_{33} - A_{33} x - \theta_1 &  B_{34} - A_{34} x- \theta_2 & B_{35} - A_{35} x \\
                    &  B_{44} - A_{44} x- \theta_3 & B_{45} - A_{45} x \\
                    &                    & B_{55} - A_{55} x \\
  \end{bmatrix} \in \Psd[3]
  \\
   S_2 &= \begin{bmatrix}
  B_{22} - A_{22} x- \theta_4 &  B_{23} - A_{23} x - \theta_5 & B_{24} - A_{24} x \\
                    &  \theta_1 - \theta_6 &  \theta_2 \\
                    &                    &  \theta_3\\
  \end{bmatrix} \in \Psd[3]
  \\
   S_1 &= \begin{bmatrix}
  B_{11} - A_{11} x &  B_{12} - A_{12} x & B_{13} - A_{13} x \\
                    &  \theta_4 & \theta_5 \\
                    &                    & \theta_6\\
  \end{bmatrix} \in \Psd[3]
\end{align*}
Notice how the overlap variables drop out if the original matrix $S$ is reassembled according to~Theorem~\ref{thm:agler} by adding the entries of each subblock.


\section{Clique Merging}
\label{sec:clique_merging}
Given an initial decomposition with (chordally completed) edge set $E$ and a set of cliques $\{\clique{1},\ldots,\clique{p}\}$, we are free to re-merge any number of cliques back into larger blocks.
This is equivalent to treating \emph{structural} zeros in the problem as \emph{numerical} zeros, leading to additional edges in the graph.
Looking at the decomposed problem in~\eqref{eq:primal_decomp_sdp} and~\eqref{eq:dual_decomp_sdp}, the effects of merging two cliques $\clique{i}$ and $\clique{j}$ are twofold:
\begin{enumerate}[1.]
  \item We replace two positive semidefinite matrix constraints of dimensions $\abs{\clique{i}}$ and $\abs{\clique{j}}$ with one constraint on a larger clique with dimension $\left| \clique{i} \cup \clique{j} \right|$, where the increase in dimension depends on the size of the overlap.
  \item We remove consistency constraints for the overlapping entries between $\clique{i}$ and $\clique{j}$, thus reducing the size of the linear system of equality constraints.
\end{enumerate}
When merging cliques these two factors have to be balanced, and the optimal merging strategy depends on the particular SDP solution algorithm employed.
In~\cite{Nakata_2003} and~\cite{Sun_2014} a clique tree is used to search for favourable merge candidates;
We will refer to their two approaches as \pkg{SparseCoLO} and the \emph{parent-child} strategy, respectively, in the following sections.   We will then propose a new merging method in Section \ref{ssec:new_cg_strategy} whose performance is superior to both methods when used in ADMM.
Given a merging strategy, Algorithm~\ref{alg:merge_function} describes how to merge a set of cliques within the set $\mc{B}$ and update the edge set $\mc{E}$ accordingly.
\begin{figure}[htb]
 \removelatexerror
\begin{algorithm}[H]
\SetKwInOut{Input}{Input}
\SetKwInOut{Output}{Output}
\SetKwRepeat{Do}{Do}{while}
\Input{A set of cliques $\mc{B}$ with edge set $\mc{E}$, a subset of cliques $\mc{B}_m = \{ \clique{m, 1} , \clique{m, 2},  \ldots, \clique{m, r}\}\subseteq \mc{B} $ to be merged.}
\Output{A reduced set of cliques $\hat{\mc{B}}$ with edge set $\hat{\mc{E}}$ and the merged clique $\clique{m}$.}
 $\hat{\mc{E}} \gets \mc{E} $\;
$\clique{m} \gets  \clique{m, 1} \cup   \clique{m, 2} \cup\cdots \cup \mc{C}_{m, r}$\;
$\hat{\mc{B}} \gets (\mc{B} \setminus \mc{B}_m) \cup \{\clique{m}\} $\;
  Remove edges  $ \{(\clique{i}, \clique{j}) \mid i \neq j, \; \clique{i}, \clique{j} \in \mc{B}_m \}$ in $\hat{\mc{E}}$\;
  Replace edges $\{(\clique{i}, \clique{j}) \mid \clique{i} \in \mc{B}_m,  \clique{j} \notin  \mc{B}_m \} $ with $(\clique{m}, \clique{j})$ in $\hat{\mc{E}}$\;
 \caption{Function $\mathrm{mergeCliques}(\mc{B}, \mc{E}, \mc{B}_m)$.}
   \label{alg:merge_function}
\end{algorithm}
\end{figure}
\subsection{Existing clique tree-based strategies}
The parent-child strategy described in~\cite{Sun_2014} traverses the clique tree in a depth-first order and merges a clique $\clique{\ell}$ with its parent clique \mbox{$\clique{\mathrm{par}(\ell)} \eqdef \mathrm{par}(\clique{\ell})$} if at least one of the two following conditions are met:
\begin{align}
   \left( \abs{\mc{C}_{\mathrm{par}(\ell)}} - \abs{\eta_\ell} \right) \left( \abs{\clique{\ell}} - \abs{\eta_\ell} \right) &\leq t_\mathrm{fill},\\
   \max \left\{ \abs{\nu_\ell}, \abs{\nu_{\mathrm{par}(\ell)}} \right\} &\leq t_\mathrm{size},
\end{align}
with heuristic parameters $t_{\mathrm{fill}}$ and $t_{\mathrm{size}}$.
These conditions keep the amount of extra fill-in and the supernode cardinalities below the specified thresholds.
The \pkg{SparseCoLO} strategy described in~\cite{Nakata_2003} and~\cite{Fujisawa_2006} considers parent-child as well as sibling relationships.
Given a parameter $\sigma > 0$, two cliques $\clique{i}, \clique{j}$ are merged if the following merge criterion holds
\begin{equation}
\label{eq:sparsecolo_criterion}
  \min \left\{ \frac{\abs{\clique{i} \cap \clique{j}}}{\abs{\clique{i}}}, \frac{\abs{\clique{i} \cap \clique{j}}}{\abs{\clique{j}}}  \right\} \geq \sigma.
\end{equation}
This approach traverses the clique tree depth-first, performing the following steps for each clique $\clique{\ell}$:
\begin{enumerate}[1.]
\item
For each clique pair \mbox{$\left\{ (\clique{i}, \clique{j}) \mid \clique{i}, \clique{j} \in \mathrm{ch}\left(\clique{\ell} \right) \right\}$}, check if~\eqref{eq:sparsecolo_criterion} holds, then:
\begin{itemize}
   \item $\clique{i}$ and $\clique{j}$ are merged, or
   \item if $(\clique{i} \cup \clique{j}) \supseteq \clique{\ell}$, then $\clique{i}$, $\clique{j}$, and $\clique{\ell}$ are merged.
 \end{itemize}
\item
For each clique pair  $\left\{ \left(\clique{i}, \clique{\ell}\right) \mid \clique{i} \in \mathrm{ch}\left(\clique{\ell}\right) \right\}$, merge $\clique{i}$ and $\clique{\ell}$  if~\eqref{eq:sparsecolo_criterion} is satisfied.
\end{enumerate}
We note that the open-source implementation of the \pkg{SparseCoLO} algorithm described in~\cite{Nakata_2003} follows the algorithm outlined here, but also employs a few additional heuristics.

\begin{figure}[htb]
  \centering
\begin{tabular}[b]{@{}cc@{}}
\begin{minipage}{.4\linewidth}
\tikzstyle{sbox}=[rectangle split,rectangle split parts=2,draw,text centered]
\centering
 \includegraphics[width=.6\textwidth]{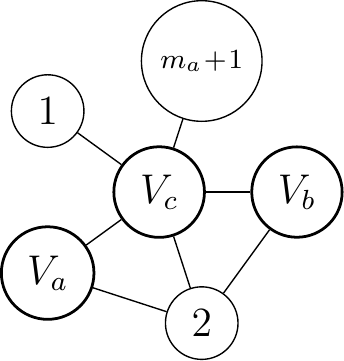}
\end{minipage}%
&
\begin{minipage}{.4\linewidth}
\centering
 \includegraphics[width=.8\textwidth]{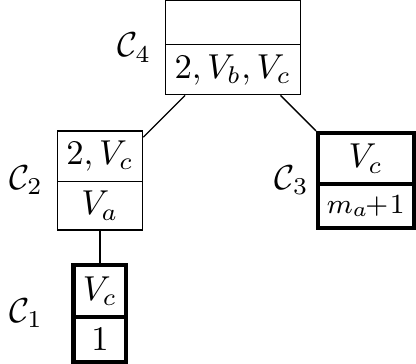}

\end{minipage}
\\
(a) & (b)
\end{tabular}

  \caption{Sparsity graph (a) that can lead to clique tree (b) with an advantageous \textquote{nephew-uncle} merge between $\clique{1}$ and $\clique{3}$.}
  \label{fig:nephew_merge}
\end{figure}

An advantage of these two approaches is that the clique tree can be computed easily and the conditions are inexpensive to evaluate.
However, a disadvantage is that choosing parameters that work well on a variety of problems and solver algorithms is difficult.
Secondly, clique trees are not unique and in some cases it is beneficial to merge cliques that are not directly related on the particular clique tree that was computed.
To see this, consider a chordal graph $G(V,E)$ consisting of three connected subgraphs:
\begin{align*}
G_a(V_a, E_a),  \text{ with } V_a &= \{3, 4, \ldots, m_a\}, \\
G_b(V_b, E_b),  \text{ with } V_b &= \{m_{a}+2, m_{a}+3, \ldots, m_b\}, \\
G_c(V_c, E_c),  \text{ with } V_c &= \{m_{b}+1, m_{b}+2, \ldots, m_c\},
\end{align*}
and some additional vertices $\{1,2,m_a + 1\}$.
The graph is connected as shown in Figure~\ref{fig:nephew_merge}(a), where the complete subgraphs are represented as nodes $V_a,V_b,V_c$.
A corresponding clique tree is shown in Figure~\ref{fig:nephew_merge}(b).
By choosing the cardinality $\abs{V_c}$, the overlap between cliques $\clique{1} = \{1, 2\} \cup V_c$ and \mbox{$\clique{3} = \{m_a+1\} \cup V_c $} can be made arbitrarily large while $\abs{V_a}$, $\abs{V_b}$ can be chosen so that any other merge is disadvantageous.
However, neither the parent-child strategy nor \pkg{SparseCoLO} would consider merging $\clique{1}$ and $\clique{3}$ since they are in a \textquote{nephew-uncle} relationship.
In fact for the particular sparsity graph in Figure~\ref{fig:nephew_merge}(a) eight different clique trees can be computed.
Only in half of them do the cliques $\clique{1}$ and $\clique{3}$ appear in a parent-child relationship.
Therefore, a merge strategy that only considers parent-child relationships would miss this favorable merge in half the cases.

\subsection{A new clique graph-based strategy}
\label{ssec:new_cg_strategy}
To overcome the limitations of existing strategies we propose a merging strategy based on the \emph{reduced clique graph} $\mc{G}(\mc{B}, \xi)$, which is defined as the union of all possible clique trees of $G$; see~\cite{Habib_2012} for a detailed discussion.
The set of vertices of this graph is given by the maximal cliques of the sparsity pattern.
We then create the edge set $\xi$ by introducing edges between pairs of cliques $(\clique{i}, \clique{j})$ if they form a separating pair $\mathcal{P}_{ij}$.
We remark that $\xi$ is a subset of the edges present in the \emph{clique intersection graph} which is obtained by introducing edges for every two cliques that intersect.
However, the reduced clique graph has the property that it remains a valid reduced clique graph of the altered sparsity pattern after performing a permissible merge between two cliques.
This is not always the case for the clique intersection graph.
For convenience, we will refer to the reduced clique graph simply as the clique graph in the following sections.
Based on the permissibility condition for edge reduction in~\cite{Habib_2012} we define a permissibility condition for clique merges:
\begin{definition}[Permissible merge]~\\
Given a reduced clique graph $\mc{G}(\mc{B}, \xi)$, a merge between two cliques $(\clique{i}, \clique{j}) \in \xi$ is permissible if for every common neighbour $\clique{k}$ it holds that $\clique{i} \cap \clique{k} =  \clique{j} \cap \clique{k}$.
\end{definition}

We further define a monotone \emph{edge weighting function} $e \colon 2^{V} \times 2^{V} \rightarrow \Re$  that assigns a weight $w_{ij}$ to each edge  $(\clique{i}, \clique{j}) \in \xi$:
\[
e\left(\clique{i}, \clique{j}\right) = w_{ij}.
\]
This function is used to estimate the per-iteration computational savings of merging a pair of cliques depending on the targeted algorithm and hardware.
It evaluates to a positive number if a merge reduces the per-iteration time and to a negative number otherwise.
For a first-order method, whose per-iteration cost is dominated by an eigenvalue factorisation with complexity $\mc{O}\bigl(\abs{\clique{}}^3\bigr)$, a simple choice would be:
\begin{equation}\label{eq:edge_weighting}
e(\clique{i}, \clique{j}) = \abs{\clique{i}}^3 + \abs{\clique{j}}^3 - \abs{\clique{i} \cup \clique{j}}^3.
\end{equation}
More sophisticated weighting functions can be determined empirically; see~Section~\ref{ssec:merging_benchmarks}.
After a weight has been computed for each edge $(\clique{i}, \clique{j})$ in the clique graph, we merge cliques as outlined in Algorithm~\ref{alg:cg_strategy}.
\begin{figure}[htb]
 \removelatexerror
\begin{algorithm}[H]
\SetKwInOut{Input}{Input}
\SetKwInOut{Output}{Output}
\SetKwRepeat{Do}{Do}{while}
\Input{A weighted clique graph $\mc{G}(\mc{B}, \xi)$. }
\Output{A merged clique graph $\mc{G}(\hat{\mc{B}}, \hat{\xi})$. }
$\hat{\mc{B}} \gets \mc{B}$ and $\hat{\xi} \gets \xi $\;
 STOP $\gets$ false\;
\While{!STOP}{
  choose permissible edge $\left(\clique{i}, \clique{j}\right)$ with maximum $w_{ij}$\;
  \eIf{$w_{ij} > 0$}{
    $\mc{B}_m \gets \{\clique{i}, \clique{j} \}$\;
    $\hat{\mc{B}}, \hat{\xi}, \clique{m} \gets \mathrm{mergeCliques}\left(\hat{\mc{B}}, \hat{\xi}, \mc{B}_m \right)$\;
    \For{each edge $(\clique{m}, \clique{\ell})   \in \hat{\xi}$}{
      update $w_{m \ell} \gets e(\clique{m}, \clique{\ell})$\;
    }
  }
  {
  STOP $\gets$ true\;
  }
}
 \caption{Clique graph-based merging strategy.}
   \label{alg:cg_strategy}
\end{algorithm}
\end{figure}
This strategy considers the edges in terms of their weights, starting with the permissible clique pair $(\clique{i}, \clique{j})$ with the highest weight $w_{ij}$.
If the weight is positive, the two cliques are merged and the edge weights for all edges connected to the merged clique $\clique{m} = \clique{i} \cup \clique{j}$ are updated.
This process continues until no edges with positive weights remain.

\begin{figure}[htb]
  \centering
\begin{tabular}[b]{@{}cc@{}}
\begin{minipage}{0.4\linewidth}
\centering
 \includegraphics[width=0.85\linewidth]{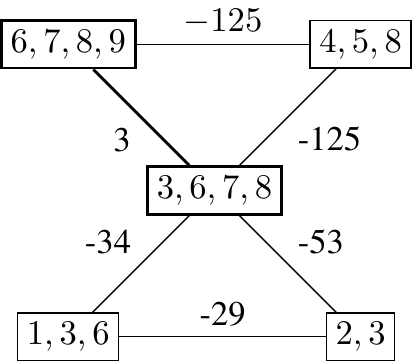}

\end{minipage}
&
\begin{minipage}{0.4\linewidth}
\centering
 \includegraphics[width=0.9\linewidth]{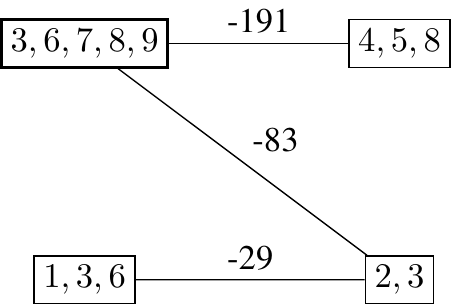}

\end{minipage}
\\
(a) & (b)
\end{tabular}
  \caption{(a)~Clique graph $\mc{G}(\mc{B}, \xi)$ of the clique tree in Figure~\ref{fig:example_structures}(c) with edge weighting function $e(\clique{i}, \clique{j}) = \abs{\clique{i}}^3 + \abs{\clique{j}}^3 - \abs{\clique{i} \cup \clique{j}}^3$ and (b)~clique graph  $\mc{G}(\hat{\mc{B}}, \hat{\xi})$ after merging the cliques $\{3, 6, 7, 8\}$ and $\{6, 7, 8, 9\}$ and updating edge weights. }
  \label{fig:merged_example}
\end{figure}

The clique graph for the clique tree in Figure~\ref{fig:example_structures}(c) is shown in Figure~\ref{fig:merged_example}(a) with the edge weighting function in~\eqref{eq:edge_weighting}.
Following Algorithm~\ref{alg:cg_strategy} the edge with the largest weight is considered first and the corresponding cliques are merged, i.e.\ $\{3, 6, 7, 8\}$ and $\{6, 7, 8, 9\}$.
Note that the merge is permissible because both cliques intersect with the only common neighbour $\{4,5,8\}$ in the same way.
The revised clique graph $\mc{G}(\hat{\mc{B}}, \hat{\xi})$ is shown in Figure~\ref{fig:merged_example}(b).
Since no edges with positive weights remain, the algorithm stops.

After~Algorithm~\ref{alg:cg_strategy} has terminated, it is possible to recompute a valid clique tree from the revised clique graph.
This can be done in two steps.
First, the edge weights in $\mc{G}(\hat{\mc{B}}, \hat{\xi})$ are replaced with the cardinality of their intersection:
\begin{equation*}
  \tilde{w}_{ij} = \abs{\clique{i} \cap \clique{j}}, \; \text{for all } (\clique{i},\clique{j}) \in \hat{\xi}.
\end{equation*}
Second, a clique tree is then given by any \emph{maximum weight spanning tree} of the newly weighted clique graph, \eg\ using Kruskal's algorithm described in~\cite{Kruskal_1956}.

Our merging strategy has some clear advantages over competing approaches.
Since the clique graph covers a wider range of merge candidates, it will consider edges that do not appear in clique tree-based approaches such as the \textquote{nephew-uncle} example in Figure~\ref{fig:nephew_merge}.
Moreover, the edge weighting function allows one to make a merge decision based on the particular solver algorithm and hardware used.
One downside is that this approach is more computationally expensive than the other methods.
However, our numerical experiments show that the time spent on finding the clique graph, merging the cliques, and recomputing the clique tree represent only a very small fraction of the total computational savings relative to other merging methods when solving SDPs.


\section{Open-Source Implementation}
\label{sec:implementation}
We have implemented our algorithm in the Conic Operator Splitting Method (\pkg{COSMO}), an open-source package written in Julia~\cite{Bezanson_2017}.
Julia allows the solver to be written in a flexible, modular and extensible way, while still maintaining the benefits of a fast compiled language.
The source code and documentation are available at
\begin{center}
\url{https://github.com/oxfordcontrol/COSMO.jl}.
\end{center}
\pkg{COSMO} offers the user two interfaces to describe the constraints of the optimisation problem: a direct interface, and an interface to the modelling languages \pkg{JuMP}~\cite{Dunning_2017} and \pkg{Convex.jl}~\cite{Udell_2014}.
These interfaces connect the solver to the Julia optimisation ecosystem which provide flexible problem description and automatic problem reformulation.

As shown in~Algorithm~\ref{alg:admm} the two main steps of the algorithm are solving a linear system and projecting onto a Cartesian product of cones.
The implementation of these two parts allows customisation by the user.
For the solution of the linear system in~\eqref{eq:sys} the user can either use the \pkg{QDLDL}~\cite{Stellato_2018} solver provided with \pkg{COSMO}, the standard sparse solver from \pkg{SuiteSparse}~\cite{Davis_2015}, or choose one of the provided interfaces to direct and indirect solvers, \eg\ \pkg{Pardiso}~\cite{Schenk_2010,Kalinkin_2015}, conjugate gradient, minimal residual method, or link their own implementation.

The second important part of the algorithm is the projection step onto a Cartesian product of convex sets.
By default \pkg{COSMO} supports the zero cone, the nonnegative orthant, the hyperbox, the second-{}order cone, the PSD cone, the exponential cone and its dual, and the power cone and its dual.
Our Julia implementation also allows the user to define their own convex cones\footnote{To allow infeasibility detection the user has to either define a convex cone, a convex compact set or a composition of the two.} and custom projection functions.
To implement a custom cone $\mc{K}_c$ the user has to provide:
\begin{itemize}
  \item a projection function that projects an input vector onto the cone
  \item a function that determines if a vector is inside the dual cone $\mathcal{K}_c^*$
  \item a function that determines if a vector is inside the recession cone of $-\mathcal{K}_c$
\end{itemize}
The latter two functions are required for our solver to implement checks for infeasibility as described in Sections~\ref{ssec:infeasibiltyCerts} and \ref{ssec:algorithmConvergence}.
An example that shows the advantages of defining a custom cone is provided in~Section~\ref{ssec:custom_set}.

The authors in~\cite{Rontsis_2019} used \pkg{COSMO}'s algorithm with a specialized implementation of the projection function for positive semidefinite constraints.
The projection method used approximate matrix eigendecompositions to significantly reduce the projection time, while maintaining all the features of \pkg{COSMO} such as scaling, infeasibility detection and interfaces to linear system solvers.
It was demonstrated that this can provide a significant, up to 20x, reduction in solve time.\\

Given a problem with multiple constraints, the projection step~\eqref{eq:proj_step} can be carried out in parallel.
This is particularly advantageous when used in combination with chordal decomposition, which typically yields a large number of smaller PSD constraints.
For the eigendecomposition involved in the projection step of a PSD constraint, the LAPACK~\cite{Anderson_1999} function \texttt{sveyr} is used, which can also utilise multiple threads.
Consequently, this leads to two-level parallelism in the computation, i.e\ on the higher level the projection functions are carried out in parallel and each projection function independently calls \texttt{sveyr}.
Determining the optimal allocation of the CPU cores to each of these tasks depends on the number of PSD constraints and their dimensions and is a difficult problem.
For the problem sets considered in section~\ref{sec:results} we achieved the best performance by running \texttt{sveyr} single-threaded and using all physical CPUs to carry out the projection functions in parallel.

Moreover, Julia's type abstraction features are used to enable the solver to solve problems of arbitrary floating-point precision.
This allows for example to reduce the memory usage of the solver for very large problems by switching to 32-bit single-precision floating-point format.

\section{Numerical Results}
\label{sec:results}

This section presents benchmark results of \pkg{COSMO} against the interior-point solver \pkg{MOSEK} v$9.0$ and the accelerated first-order ADMM solver \pkg{SCS} v$2.1.1$.
When applied to a quadratic program, \pkg{COSMO}'s main algorithm becomes very similar to the first-order QP solver~\pkg{OSQP}.
To test the performance penalty of using a pure Julia implementation against a similar C implementation we also compare our solver against~\pkg{OSQP} v$0.6.0$ on QP problems.

We selected a number of problem sets to test different aspects of \pkg{COSMO}.
The advantage of supporting a quadratic cost function in a conic solver is shown by solving QPs from the Maros and M\'{e}sz\'{a}ros QP repository~\cite{Maros_1999} in Section~\ref{ssec:maros} and SDPs with quadratic objectives in the form of nearest correlation matrix problems in Section~\ref{ssec:nearest_correlation_matrix}.

To highlight the advantages of implementing custom constraints, we consider a problem set with doubly-stochastic matrices in Section~\ref{ssec:custom_set}.
We then show how chordal decomposition can speed up the solver for SDPs that exhibit a block-arrow sparsity pattern in Section~\ref{ssec:block_arrow_SDPs}.

The performance of chordal decomposition is further explored in Section~\ref{ssec:merging_benchmarks} by solving large structured problems from the SDPLib benchmark set~\cite{Borchers_1999} as well as some non-chordal SDPs generated with sparsity patterns from the \pkg{SuiteSparse} Matrix Collections~\cite{Davis_2011}.
Using the same problems we additionally evaluate the performance of different clique merging strategies.

All the experiments were carried out on a computing node of the University of Oxford ARC-HTC cluster with 16 logical Intel Xeon E5-2560 cores and \SI{64}{GB} of DDR3 RAM.
All the problems were run using Julia v$1.3$ and the problems were passed to the solvers via MathOptInterface~\cite{Legat_2020}.

To evaluate the accuracy of the returned solution we compute three errors adapted from the \mbox{DIMACS} error measures for SDPs~\cite{Johnson_2000}:
\begin{equation}
\epsilon_1 = \frac{\norm{A_ax - b_a}_2}{1 + \norm{b_a}_2}, \quad \epsilon_2 = \frac{\norm{Px + q - A_a^\top y_a}_2}{1 + \norm{q}_2}, \quad \epsilon_3 = \frac{ \abs{x^\top P x +   q^\top x - b_a^\top y_a} }{ 1 + \abs{q^\top x} + \abs{b_a^\top y_a}},
\end{equation}
where $A_a$, $b_a$ and $y_a$ correspond to the rows of $A$, $b$ and $y$ that represent active constraints.
This is to ensure meaningful values even if the problem contains inactive constraints with very large, or possibly infinite, values $b_i$.
The maximum of the three errors for each problem and solver is reported in the results below.

We configured \pkg{COSMO}, \pkg{MOSEK}, \pkg{SCS} and \pkg{OSQP} to achieve an accuracy of $\epsilon = \num{e-3}$.
We set the maximum allowable solve time for the Maros and M\'{e}sz\'{a}ros problems to $\SI{5}{\minute}$ and to $\SI{30}{\minute}$ for the other problem sets.
All other solver parameters were set to the solvers' standard configurations.
\pkg{COSMO} uses a Julia implementation of the \pkg{QDLDL} solver to factor the quasi-definite linear system.
Similarly, we configured \pkg{SCS} to use a direct solver for the linear system.

\subsection{Maros and M\'{e}sz\'{a}ros QP test set}
\label{ssec:maros}
The Maros and M\'{e}sz\'{a}ros test problem set~\cite{Maros_1999} is a repository of challenging convex QP problems that is widely used to compare the performance of QP solvers.
For comparison metrics we compute the failure rate, the number of fastest solve time and the normalized shifted geometric mean for each solver.
The shifted geometric mean is more robust against large outliers (compared to the arithmetic mean) and against small outliers (compared to the geometric mean) and is commonly used in optimisation benchmarks; see~\cite{Stellato_2018,Mittelmann}.
The shifted geometric mean $\mu_{g, s}$ is defined as:
\begin{equation}
  \mu_{g, s} \eqdef \sqrt[n]{\prod_{p}(t_{p, s} + \mathrm{sh}) - \mathrm{sh} }
\end{equation}
with total solver time $t_{p, s}$ of solver $s$ and problem $p$, shifting factor $\mathrm{sh}$ and size of the problem set $n$.
In the reported results a shifting factor of $sh = 10$ was chosen and the maximum allowable time $t_{p, s} = \SI{300}{s}$ was used if solver $s$ failed on problem $p$.
Lastly, we normalize the shifted geometric mean for solver $s$ by dividing by the geometric mean of the fastest solver.
The failure rate $f_{r, s}$ is given by the number of unsolved problems compared to the total number of problems in the problem set.
As unsolved problems we count instances where the algorithm does not converge within the allowable time or fails during the setup or solve phase.
Table~\ref{tb:maros_meszaros} shows the normalized shifted geometric mean and the failure rate for each solver. Additionally, the number of cases where solver $s$ was the fastest solver is shown.

\begin{table}[htb]
\begin{center}
\sisetup{round-mode=places
        ,round-precision=3
        ,scientific-notation=false}
\caption{Normalized shifted geometric mean and failure rates for solvers tested on the Maros and M\'{e}sz\'{a}ros QP test set.}\label{tb:maros_meszaros}
\begin{threeparttable}
   \begin{tabular}{r l l l l}
     \toprule
      & OSQP & COSMO & MOSEK & SCS\\
 \midrule
 Normalized shifted geometric mean & $\num{1.0}$ & $\num{1.1688398708116283}$ & $\num{1.896856715452622}$ & $\num{75.015173814748}$\\
  Number of fastest solve time & $\num{75}$ & $\num{27}$ & $\num{33}$ & $\num{0}$\\
  Failure rates $f_{r,s}$ [\%] & $\num{4.444444444444445}$ & $\num{5.185185185185185}$ & $\num{10.37037037037037}$ & $\num{83.7037037037037}$\\

   \bottomrule
   \end{tabular}
 \end{threeparttable}
\end{center}
\end{table}
\pkg{OSQP} shows the best performance in terms of lowest failure rates, number of fastest solves and in the shifted geometric mean of solve times.
\pkg{COSMO} follows very closely.
The shifted geometric mean of \pkg{MOSEK} seems to suffer from a higher failure rate compared to \pkg{OSQP}/\pkg{COSMO}, and \pkg{SCS} fails on a large number of problems.
The higher failure rate could be due to the necessary transformation into a second-order-cone problem.

\begin{figure}[h]
  \centering
      \includegraphics[width=\textwidth]{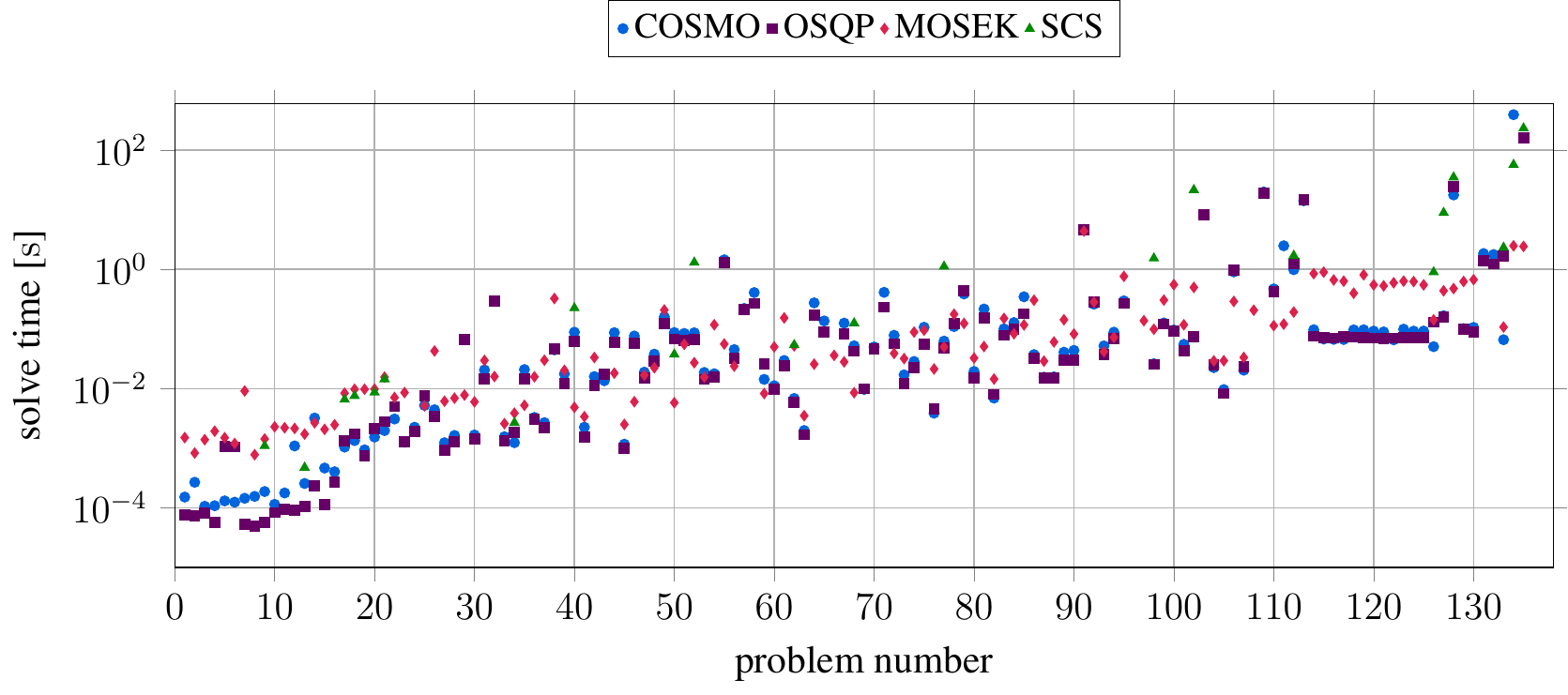}
  \caption{Solve time of benchmarked solvers for problems of the Maros and M\'{e}sz\'{a}ros QP problem set. Only problem results classified as solved are shown. The problems are ordered by increasing number of non-zeros in the constraint matrix.}
  \label{fig:bench_maros}
\end{figure}

For this problem set of QPs \pkg{COSMO}'s algorithm reduces, with some minor differences, to the algorithm of \pkg{OSQP}.
Consequently, this benchmark is useful to evaluate the performance penalty that \pkg{COSMO} pays due to its implementation in the higher-order language Julia.
The results in~Table~\ref{tb:maros_meszaros} show that the performance difference is very small.
This can also be seen by looking at the solve times of each solver for increasing problem dimension, as shown in~Figure~\ref{fig:bench_maros}.

\pkg{COSMO} and \pkg{OSQP} have very similar solve times, aside from very small problems that are solved in under \SIrange{1e-5}{1e-4}{\s}.  This difference is primarily due to overheads incurred  from features in our Julia implementation that support more than one constraint type during problem setup.  The marginally better resulting performance of OSQP for the smallest problems in the test set is the reason that OSQP is the faster solver in a larger number of cases in Table~\ref{tb:maros_meszaros}.

\subsection{Custom convex cones}
\label{ssec:custom_set}
In many cases writing a custom solver algorithm for a particular problem can be faster than using available solver packages if a particular aspect of the problem structure can be exploited to speed up parts of the computations.
As mentioned earlier, \pkg{COSMO} supports user customisation by allowing the definition of new convex cones.
This is useful if constraints of the problem can be expressed using this new convex cone and a fast projection method onto the cone exists.
A fast specialized projection method in an ADMM framework has for example been used by the authors in~\cite{Barman_2013} to solve the error-correcting code decoding problem.

To demonstrate the advantage of custom convex cones, consider the problem of finding the doubly stochastic matrix that is nearest, in the Frobenius norm, to a given symmetric matrix $C \in \Sym[n]$.
Doubly stochastic matrices are used for instance in spectral clustering~\cite{Zass_2007} and matrix balancing~\cite{Rao_2014}.
A specialized algorithm for this problem type has been recently discussed by the authors in~\cite{Rontsis_2019a}.
Doubly stochastic matrices have the property that all rows and columns each sum to one and all entries are nonnegative.
The nearest doubly stochastic matrix $X$ can be found by solving the following optimisation problem:
\begin{equation}
  \begin{array}{ll}
    \mbox{minimize}   & \textstyle{\frac{1}{2}} \norm{X- C}_F^2\\
    \mbox{subject to} & X_{ij} \geq 0\\
                      & X \mathbf{1} = \mathbf{1} \\
                      & X^\top \mathbf{1} = \mathbf{1},
  \end{array}
  \label{eq:dsm}
\end{equation}
with symmetric real  matrix $C \in \Sym[n]$ and decision variable $X \in \Re^{n\times n}$.
This problem can be solved as a QP in the following form using equality and inequality constraints:
\begin{equation}
  \begin{array}{ll}
    \mbox{minimize}   & \frac{1}{2}(x^\top x - 2c^\top x + c^\top c)\\
    \mbox{subject to} & \begin{bmatrix}
    \mathbf{1}^\top_n \otimes I_n \\
    I_n \otimes \mathbf{1}^\top_n  \\
    - I_{n^2}
    \end{bmatrix} x + s = \begin{bmatrix}
      \mathbf{1}_{2n} \\ \mathbf{1}_{2n} \\ \mathbf{0}_{n^2}
    \end{bmatrix}\\
      & s \in \{0\}^{4n} \times \Re_+^{n^2},
  \end{array}
  \label{eq:dsm_qp}
\end{equation}
with $x = \mathrm{vec}(X)$ and $c = \mathrm{vec}(C)$
However, the problem can be written in a more compact form by using a custom projection function to project the matrix iterate onto the affine set of matrices $\mathcal{C}_{\sum}$, whose rows and columns each sum to one.
In general the projection of vector $s \in \Re^n$ onto the affine set $\mc{C}_a = \{ s \in \Re^n \mid As = b \}$ is given by:
\begin{equation}
  \Pi_{\mc{C}_a}(s) = s - A^\top \left(A A^\top \right)^{-1} (A s - b),
\end{equation}
where $A$ is assumed to have full rank.
In the case of $\mc{C}_a = \mc{C}_{\sum}$ we can exploit the fact that the inverse of $AA^\top$ can be efficiently computed.
The projection can be carried out as described in~Algorithm~\ref{alg:proj_dsm}; see  Appendix~\ref{app_project} for a derivation.
\begin{figure}[htb]
 \removelatexerror
\begin{algorithm}[H]
$A_r = \mathbf{1}_n^\top \otimes I_n$ \quad $A_c = \begin{bmatrix} I_{n-1} \otimes \mathbf{1}_n^\top & \mathbf{0}_{n-1 \times n} \end{bmatrix}$\;
$A = \begin{bmatrix} A_r & A_c \end{bmatrix}^\top$\;
$r = [r_1, r_2]^\top = A s - \mathbf{1}_{2n- 1}$\;
$\eta_2 = \frac{1}{n} \left( I_{n-1} + \mathbf{1}_{n-1}\mathbf{1}_{n-1}^\top \right) \cdot \left(  r_2 - \frac{1}{n} \mathbf{1}_{n-1} \mathbf{1}_{n}^\top r_1 \right)$\;
$\eta_1 = \frac{ 1}{n} \left( r_1 - \mathbf{1}_{n} \mathbf{1}_{n-1}^\top \eta_2 \right)$\;
$\eta = [\eta_1, \eta_2]^\top$\;
$\Pi_{\mc{C}_{\sum}}(s) = s - A^\top \eta $\;

 \caption{Projection of $s \in \Re^n$ onto $C_{\sum}$}
   \label{alg:proj_dsm}
\end{algorithm}
\end{figure}
Notice that~Algorithm~\ref{alg:proj_dsm} can be implemented efficiently without ever assembling and storing $A$ and $\mathbf{1} \mathbf{1}^\top$.
By using the custom convex set $\mathcal{C}_{\sum}$ and the corresponding projection function, \eqref{eq:dsm} can now be rewritten as:
\begin{equation}
  \begin{array}{ll}
    \mbox{minimize}   & (1/2)(x^\top x - 2c^\top x + c^\top c)\\
    \mbox{subject to} & \begin{bmatrix}
    -I_{n^2} \\
    - I_{n^2}
    \end{bmatrix} x + s = \mathbf{0}_{2 n^2} \\
      & s \in \mathcal{C}_{\sum} \times \Re_+^{n^2}.
  \end{array}
  \label{eq:dsm_custom}
\end{equation}
The sparsity pattern of the new constraint matrix $A$ only consists of two diagonals and the number of non-zeros reduces from $3n^2$ to $2n^2$.
We expect this to reduce the initial factorisation time of the linear system in~\eqref{eq:sys} as well as the forward- and back-substitution steps.

\begin{figure}[h]
  \centering
     \includegraphics[width=0.9\textwidth]{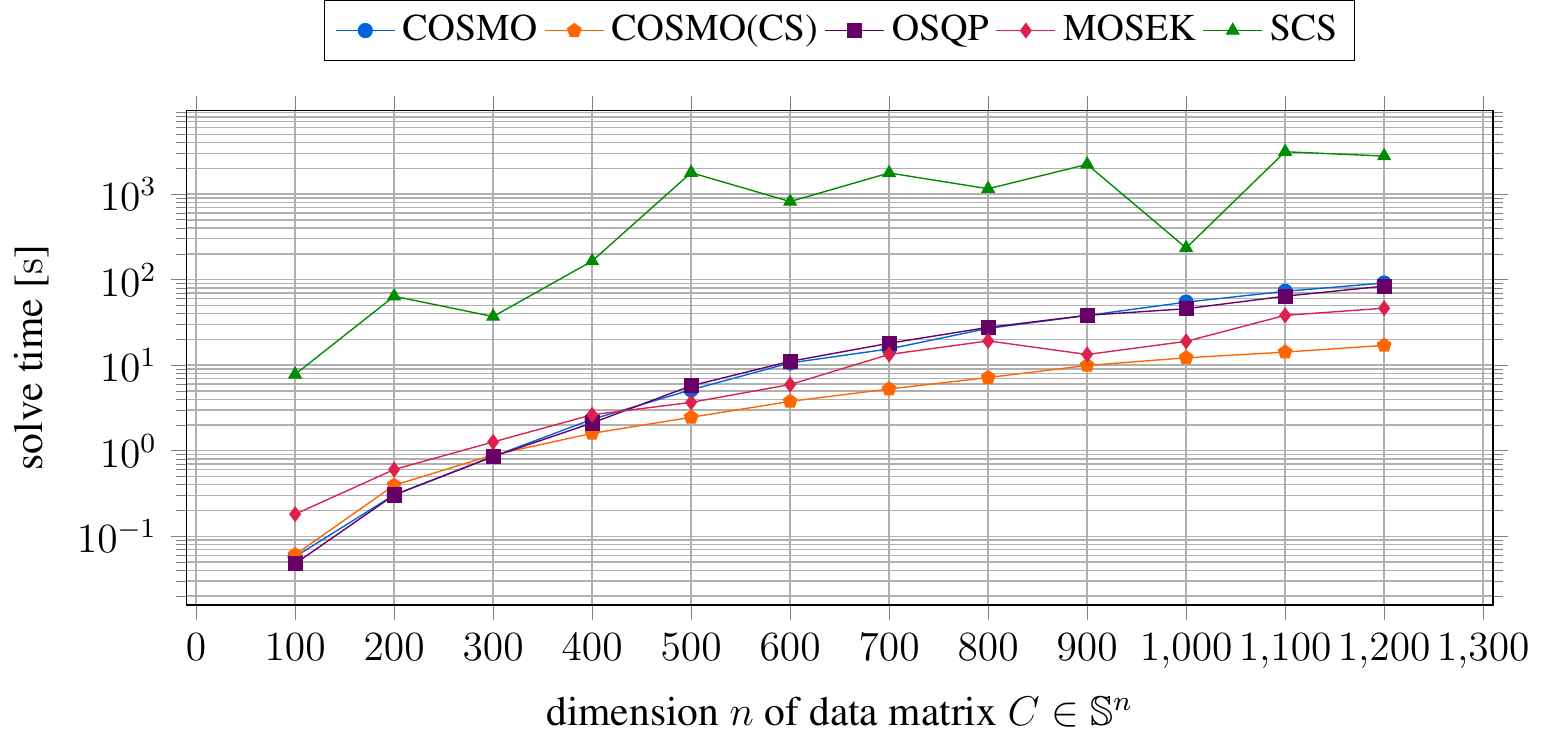}
  \caption{Solve time of benchmarked solvers for increasing problem size of doubly stochastic matrix problems. The orange line shows the solve time of \pkg{COSMO(CS)} with a custom convex set and projection function.}
  \label{fig:bench_dsm}
\end{figure}

Figure~\ref{fig:bench_dsm} shows the total solve time of all the solvers for problem~\eqref{eq:dsm} with randomly generated dense matrix $C$ with $C_{ij} \sim \mathcal{U}(0, 1)$ and increasing matrix dimension. Additionally, we show the solve time for \pkg{COSMO} in the problem form~\eqref{eq:dsm} and with a specialized custom set and projection function as in~\eqref{eq:dsm_custom}.
It is not surprising that \pkg{COSMO} and \pkg{OSQP} scale in the same way for this problem type.
\pkg{MOSEK} is slightly slower for smaller problem dimension and overtakes \pkg{COSMO}/\pkg{OSQP} for problems of dimensions $n \geq 500$.
This might be due to fact that \pkg{MOSEK} uses a faster multi-threaded linear system solver while \pkg{OSQP}/\pkg{COSMO} rely on the single-threaded solver \pkg{QDLDL}.
The longer solve time of \pkg{SCS} is due to slow convergence of the algorithm for this problem type.
Furthermore, when the problem is solved with a custom convex set as in~\eqref{eq:dsm_custom} \pkg{COSMO}(CS) is able to outperform all other solvers.
Table~\ref{tb:doubly_stochastic} shows the total solve time and the factorisation time of the two versions of \pkg{COSMO} for small, medium and large problems.
As predicted the lower solve time can be mainly attributed to the faster factorisation time.

\begin{table}[htb]
\begin{center}
\sisetup{round-mode=places
        ,round-precision=3
        ,scientific-notation=false}
\caption{Solve times and factorisation times of \pkg{COSMO} and \pkg{COSMO(CS)} for small, medium and large doubly stochastic matrix problems.}\label{tb:doubly_stochastic}
\begin{threeparttable}
   \begin{tabular}{l l l l l}
     \toprule
      & \multicolumn{2}{c}{Factorisation time (\si{s}) } & \multicolumn{2}{c}{Solve time (\si{s})} \\
          \cmidrule(lr){2-3}  \cmidrule(lr){4-5}
      $n$ & COSMO & COSMO(CS)\tnote{1}  &  COSMO & COSMO(CS)\tnote{1} \\
     \midrule
     $100$ & $\num{0.017689224}$ & $\num{0.006380161}$ & $\num{0.05763697624206543}$ & $\num{0.06084895133972168}$ \\ 
$400$ & $\num{0.994922432}$ & $\num{0.137531562}$ & $\num{2.3287200927734375}$ & $\num{1.5964269638061523}$ \\ 
$800$ & $\num{13.551001083}$ & $\num{0.552448564}$ & $\num{26.932005882263184}$ & $\num{7.163206100463867}$ \\ 
$1200$ & $\num{52.716781178}$ & $\num{1.321927381}$ & $\num{91.27830505371094}$ & $\num{17.032224893569946}$ \\ 

   \bottomrule
   \end{tabular}
 \begin{tablenotes}
 \scriptsize{%
    \item[1] solving~\eqref{eq:dsm_custom} with a custom convex set $\mathcal{C}_{\sum}$ and projection function
    }%
 \end{tablenotes}
 \end{threeparttable}
\end{center}
\end{table}

\subsection{Nearest correlation matrix}
\label{ssec:nearest_correlation_matrix}
Consider the problem of projecting a matrix $C$ onto the set of correlation matrices, i.e.\ real symmetric positive semidefinite matrices with diagonal elements equal to $1$.
This problem is for example relevant in portfolio optimisation~\cite{Higham_2002}.
The correlation matrix of a stock portfolio might lose its positive semidefiniteness due to noise and rounding errors of previous data manipulations.
Consequently, it is of interest to find the nearest correlation matrix $X$ to a given data matrix $C \in \mathbb{R}^{n\times n}$.
The problem is given by:
\begin{equation}
  \begin{array}{ll}
\mbox{minimize} & \textstyle{\frac{1}{2}} \norm{X- C}_F^2
\\
\mbox{subject to} & X_{ii} = 1,\quad i=1, \dots, n \\
 & X \in \mathbb{S}_+^n.
  \end{array}
\label{eq:corr1}
\end{equation}
In order to transform the problem into the standard form~\eqref{eq:primal} used by \pkg{COSMO}, $C$ and $X$ are vectorized and the objective function is expanded:
\begin{equation}
  \begin{array}{ll}
\mbox{minimize}   &  (1/2)(x^\top x - 2c^\top x + c^\top c)\\
\mbox{subject to} &  \begin{bmatrix}E \\ -I \end{bmatrix} x +s  = \begin{bmatrix} \mathbf{1}_{n} \\ \mathbf{0}_{n^2} \end{bmatrix}\\
& s \in \{0\}^n \times  \mathcal{S}_+^n,
  \end{array}
  \label{eq:corr2}
\end{equation}with $c=\text{vec}(C) \in \mathbb{R}^{n^2}$ and $x=\text{vec}(X) \in \mathbb{R}^{n^2}$. Here $E \in \mathbb{R}^{n\times n^2}$ is a matrix that extracts the $n$ diagonal entries $X_{ii}$ from its vectorized form $x$.

\begin{figure}[h]
  \centering
    \includegraphics[width=0.8\textwidth]{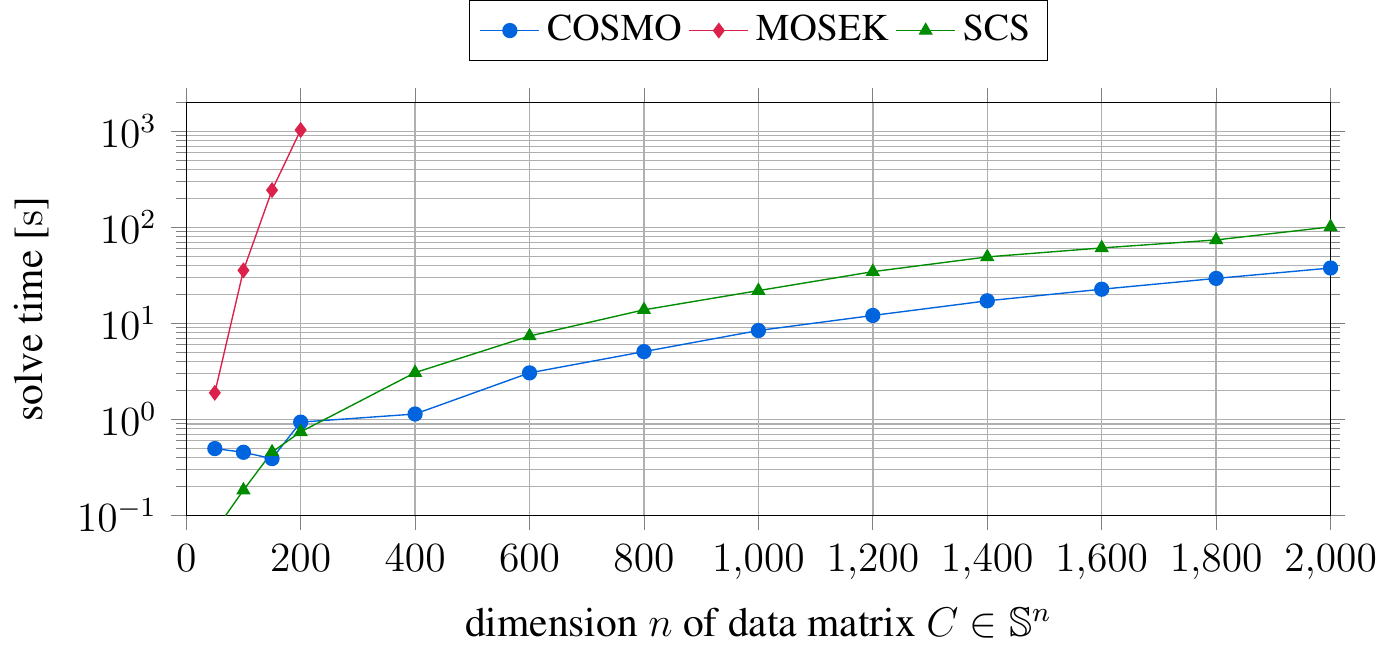}
  \caption{Solve time of benchmarked solvers for increasing problem size of nearest correlation matrix problems. The results for \pkg{MOSEK} are shown until they exceeded the time limit of $\SI{30}{min}$.}
  \label{fig:bench2}
\end{figure}

For the benchmark problems we randomly sample the data matrix $C$ with entries $C_{i,j} \sim \mc{U}(-1,1)$ from a uniform distribution.
Figure~\ref{fig:bench2} shows the benchmark results for increasing matrix dimension $n$.
Unsurprisingly, the first-order methods \pkg{SCS} and \pkg{COSMO} outperform the interior-point solver \pkg{MOSEK} for these large SDPs.
Furthermore, for larger problems the solve times of \pkg{COSMO} and \pkg{SCS} scale in a similar way.
\pkg{COSMO} seems to benefit from directly supporting the quadratic objective term in the problem while \pkg{SCS} has to transform it into an additional second-order-cone constraint. This increases the factorisation time and the projection time; see Table~\ref{tb:correlation}.

\subsection{Block-arrow sparse SDPs}
\label{ssec:block_arrow_SDPs}
\begin{figure}[h]
  \centering
   \includegraphics[width=0.3\textwidth]{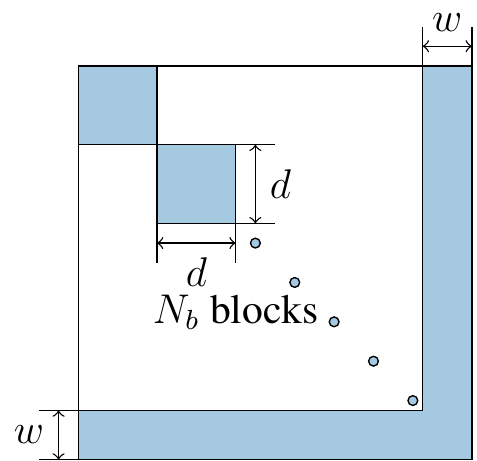}
\caption{Parameters of block-arrow sparsity pattern. The shaded area represents the non-zeros of the sparsity pattern.}
\label{fig:blockarrow}
\end{figure}

\begin{figure}[ht]
  \centering
    \includegraphics[width=0.7\textwidth]{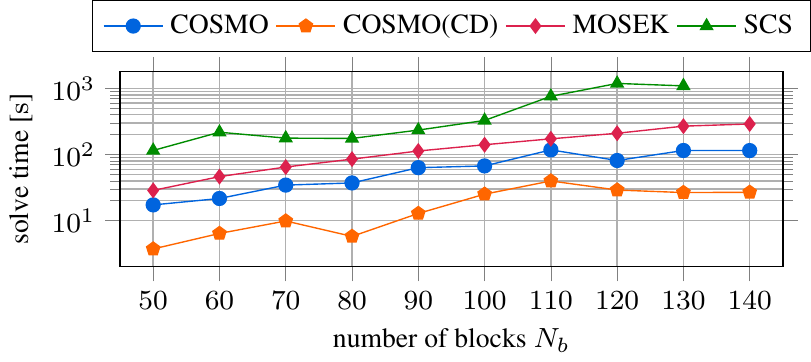}
  \caption{Solve time for increasing number of blocks $N_b$ of block-arrow sparsity pattern.}
  \label{fig:bench3_varL}
\end{figure}

\begin{figure}[ht]
  \centering
   \includegraphics[width=0.7\textwidth]{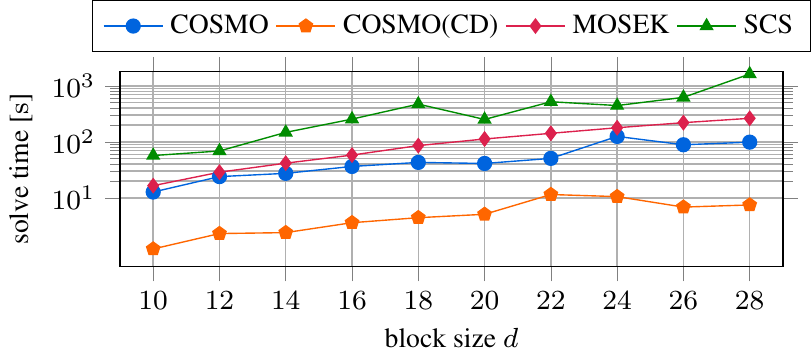}
  \caption{Solve time for increasing block size $d$ of block-arrow sparsity pattern.}
  \label{fig:bench3_varD}
\end{figure}

To demonstrate the benefits of the chordal decomposition discussed in Section~\ref{sec:chordal_decomp}, we consider randomly generated SDPs of the form~\eqref{eq:dual_sdp} with a block-arrow aggregate sparsity pattern similar to test problems in~\cite{Zheng_2019,Andersen_2010a}. Figure~\ref{fig:blockarrow} shows the sparsity pattern of the PSD constraint.
The sparsity pattern is generated based on the following parameters: block size $d$, number of blocks $N_b$ and width of the arrow head $w$.
Note that the graph corresponding to the sparsity pattern is always chordal and that, for this sparsity pattern, clique merging yields no benefit.
In the following we study the effects of independently increasing the block size $d$ and the number of blocks $N_b$. The parameters for the two test cases are:
\begin{itemize}
  \item Varying the number of blocks: $N_b=50,60,\ldots,140$, $d=10$, $w=20$, and $m=100$.
  \item Varying the block size: $d=10,12,\ldots,28$, $N_b=50$, $w=10$, and $m=100$.
\end{itemize}
The solve times for all solvers are shown in Figure~\ref{fig:bench3_varL}, Figure~\ref{fig:bench3_varD} and in~Table~\ref{tb:blockarrow}.
The line \pkg{COSMO(CD)} corresponds to the solver with chordal decomposition enabled.

The figures show that both \pkg{COSMO} with and without chordal decomposition solve this problem type consistently faster than \pkg{MOSEK} and \pkg{SCS}.
The reason why \pkg{COSMO} performs better than the other first order solver \pkg{SCS} can be explained by the significantly lower number of iterations (Table~\ref{tb:blockarrow}).
Furthermore, one can see that in both cases the solver time for each solver rises when the number of blocks and the block sizes are increased.
The increase is smaller for \pkg{COSMO(CD)} which is more affected by the number of iterations than the problem dimension.

\subsection{Non-chordal problems with clique merging}
\label{ssec:merging_benchmarks}
\begin{figure}
\centering
\begin{tabular}{ccc}
  \includegraphics[width=0.3\textwidth]{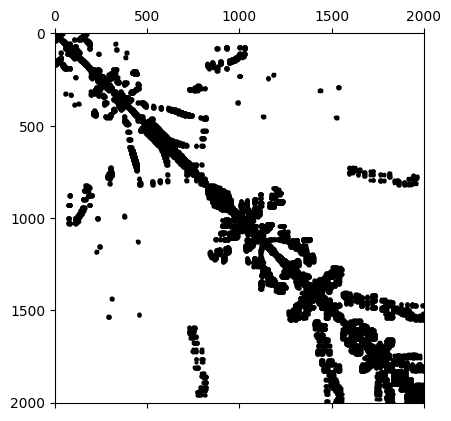} &   \includegraphics[width=0.3\textwidth]{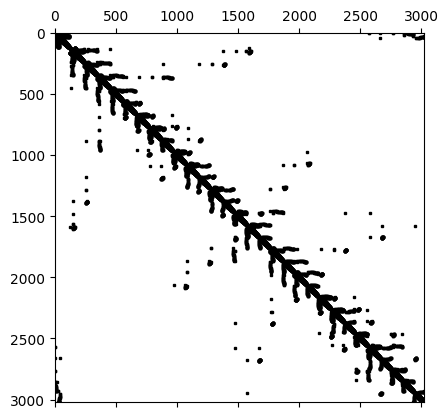} &   \includegraphics[width=0.3\textwidth]{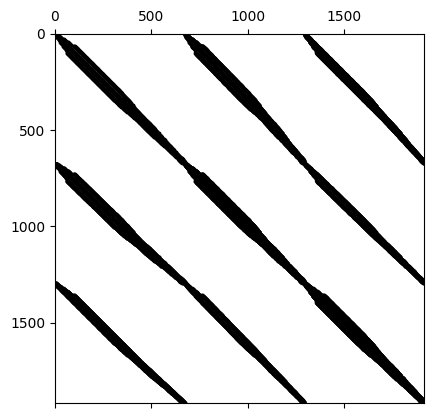}\\
(a) rs35 & (b) rs200 & (c) rs228 \\[6pt]
 \includegraphics[width=0.3\textwidth]{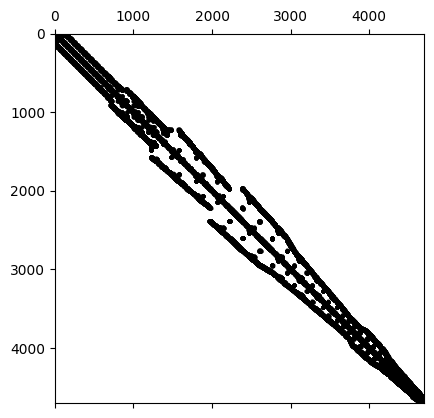} &   \includegraphics[width=0.3\textwidth]{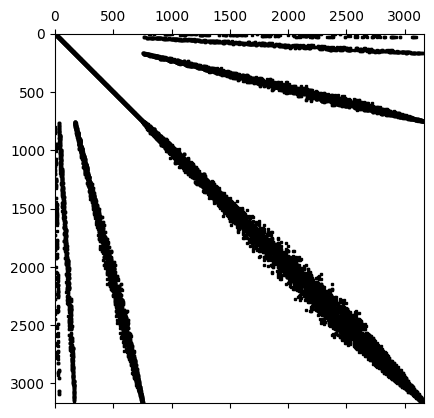} & \includegraphics[width=0.3\textwidth]{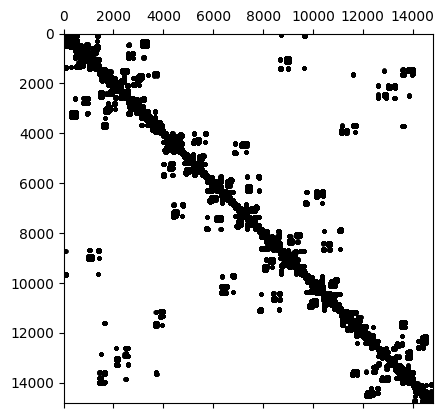}\\
(d) rs365 & (e) rs828 & (f) rs1184\\[6pt]
  \includegraphics[width=0.3\textwidth]{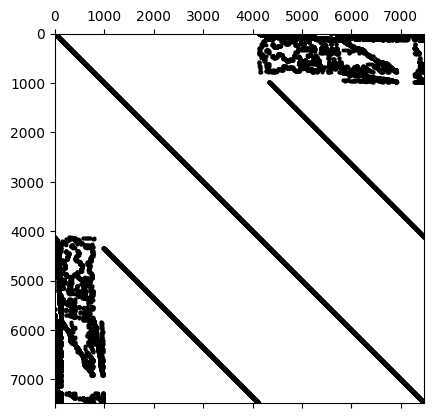} & \includegraphics[width=0.3\textwidth]{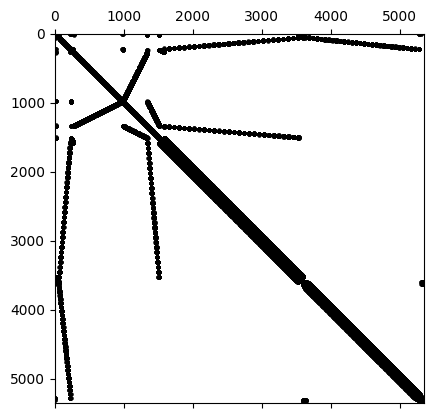}\\
 (g) rs1555 & (h) rs1907
\end{tabular}
\caption{Aggregate sparsity pattern of non-chordal SDPs created from matrices of the \pkg{SuiteSparse} Matrix Collection. The patterns are labeled with their ID number.}
\label{fig:sparsity_patterns}
\end{figure}

To compare our proposed clique graph-based merge approach with the clique tree-based strategies of~\cite{Nakata_2003} and~\cite{Sun_2014},
all three methods discussed in~Section~\ref{sec:clique_merging} were used to preprocess large sparse SDPs from SDPLib, a collection of SDP benchmark problems~\cite{Borchers_1999}.
This problem set contains maximum cut problems, SDP relaxations of quadratic programs and Lov\'{a}sz theta problems.
Moreover, we consider a set of test SDPs generated from (non-chordal) sparsity patterns of matrices from the \pkg{SuiteSparse} Matrix Collections~\cite{Davis_2011}.
The sparsity patterns for these problems are shown in Figure~\ref{fig:sparsity_patterns}.

Both problem sets were used in the past to benchmark structured SDPs~\cite{Zheng_2019,Andersen_2010a}.
This section discusses how the different decompositions affect the per-iteration computation times of the solver.
In a second step we compare the solver time of \pkg{COSMO} with our clique graph merging strategy to those of \pkg{MOSEK} and \pkg{SCS}.

\begin{figure}[htb]
\centering
\includegraphics[width=0.7\textwidth]{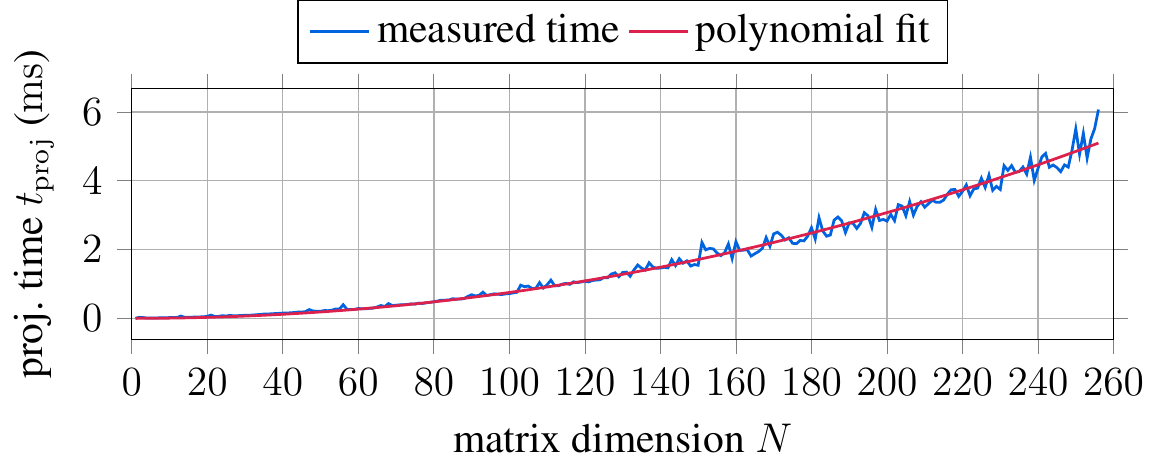}
\caption{Measured and estimated relationship between matrix size and execution time of the projection function in \pkg{COSMO}.}
  \label{fig:proj_estimate}
\end{figure}

For the strategy described in~\cite{Nakata_2003} we used the \pkg{SparseCoLO} package to decompose the problem.
The parent-child method discussed in~\cite{Sun_2014} and the clique graph based method described in~Section~\ref{ssec:new_cg_strategy} are available as options in \pkg{COSMO}.
We further investigate the effect of using different edge weighting functions.
The major operation affecting the per-iteration time is the projection step.
This step involves an eigenvalue decomposition of the matrices corresponding to the cliques.
Since the eigenvalue decomposition of a symmetric matrix of dimension $N$ has a complexity of $\mc{O}\left( N^3 \right)$, we define a \emph{nominal} edge weighting function as in~\eqref{eq:edge_weighting}.
However, the exact relationship will be different because the projection function involves copying of data and is affected by hardware properties such as cache size.
We therefore also consider an empirically \emph{estimated} edge weighting function.
To determine the relationship between matrix size and projection time, the execution time of the relevant function inside \pkg{COSMO} was measured for different matrix sizes.
We then approximated the relationship between projection time, $t_{\mathrm{proj}}$, and matrix size, $N$, as a polynomial:
\[
  t_{\mathrm{proj}}(N) = aN^3 + bN^2,
\]
where $a,b$ were estimated using least squares (Figure~\ref{fig:proj_estimate}).
The estimated weighting function is then defined as
\begin{equation}
\label{eq:estimated_edge_weighting}
  e(\clique{i}, \clique{j}) = t_{\mathrm{proj}}\left( \abs{\clique{i}} \right) + t_{\mathrm{proj}} \left( \abs{\clique{j}} \right)  - t_{\mathrm{proj}} \left( \abs{\clique{i} \cup \clique{j}} \right).
\end{equation}

Six different cases were considered: no decomposition (\texttt{NoDe}), no clique merging (\texttt{NoMer}), decomposition using \pkg{SparseCoLO} (\texttt{SpCo}), parent-child merging (\texttt{ParCh}), and the clique graph-based method with nominal edge weighting (\texttt{CG1}) and estimated edge weighting (\texttt{CG2}).
\pkg{SparseCoLO} was used with default parameters.
All cases were run single-threaded.
Since the per-iteration projection times for some problems lie in the millisecond range every problem was benchmarked ten times and the median values are reported.
Table~\ref{tb:benchmark_results} shows the solve time, the mean projection time, the number of iterations, the number of cliques after merging, and the maximum clique size of the sparsity pattern for each problem and strategy.
The solve time includes the time spent on decomposition and clique merging.
We do not report the total solver time when \pkg{SparseCoLO} was used for the decomposition because this has to be done in a separate preprocessing step in \pkg{MATLAB} which was orders of magnitude slower than the other methods.


\begin{table}[htb]
\begin{center}
   \sisetup{round-mode=places
           ,round-precision=2
           ,scientific-notation=true}
\caption{Benchmark results for non-chordal sparse SDPs from SDPLib and SDPs generated with sparsity patterns from the \pkg{SuiteSparse} Matrix Collection for different merging strategies.}\label{tb:benchmark_results}
\begin{threeparttable}
\resizebox{\textwidth}{!}{%
   \begin{tabular}{l l l l l l l l l l l l l}
     \toprule
     problem & \multicolumn{6}{c}{Median Solve time (\si{s}) } & \multicolumn{6}{c}{Median Projection time (\si{ms})}\\
         \cmidrule(lr){2-7}\cmidrule(lr){8-13}
      & NoDe\tnote{1} & NoMer\tnote{2} & SpCo\tnote{3} & ParCh\tnote{4} & CG1\tnote{5} & CG2\tnote{6} & NoDe & NoMer & SpCo & ParCh & CG1 & CG2 \\
     \midrule
     maxG11 & $35.85$ & $4.69$ & $***$ & $3.33$ & $3.09$ & \cellcolor{\myCellColor} $\mathbf{2.93}$ & $143.3$ & $18.4$ & $12.7$ & \cellcolor{\myCellColor} $\mathbf{10.5}$ & $13.4$ & $11.2$ \\ 
maxG32 & $407.82$ & $25.95$ & $***$ & $15.62$ & \cellcolor{\myCellColor} $\mathbf{14.38}$ & $19.56$ & $1257.7$ & $73.2$ & $73.1$ & $49.9$ & $51.6$ & \cellcolor{\myCellColor} $\mathbf{43.1}$ \\ 
maxG51 & $40.61$ & $32.8$ & $***$ & $9.56$ & \cellcolor{\myCellColor} $\mathbf{6.93}$ & $9.62$ & $245.8$ & $230.4$ & $219.3$ & $122.2$ & $65.5$ & \cellcolor{\myCellColor} $\mathbf{58.5}$ \\ 
mcp500-1 & $21.46$ & $1.33$ & $***$ & $0.7$ & $0.89$ & \cellcolor{\myCellColor} $\mathbf{0.6}$ & $56.0$ & $7.7$ & $7.4$ & \cellcolor{\myCellColor} $\mathbf{4.6}$ & $6.1$ & $4.9$ \\ 
mcp500-2 & $13.55$ & $12.28$ & $***$ & $8.37$ & $2.6$ & \cellcolor{\myCellColor} $\mathbf{2.42}$ & $54.9$ & $45.4$ & $32.1$ & $28.0$ & $14.3$ & \cellcolor{\myCellColor} $\mathbf{11.4}$ \\ 
mcp500-3 & $11.28$ & $32.45$ & $***$ & $30.13$ & $7.01$ & \cellcolor{\myCellColor} $\mathbf{5.14}$ & $51.1$ & $104.9$ & $105.0$ & $83.4$ & $28.2$ & \cellcolor{\myCellColor} $\mathbf{21.2}$ \\ 
mcp500-4 & $14.57$ & $72.27$ & $***$ & $13.87$ & \cellcolor{\myCellColor} $\mathbf{5.84}$ & $7.67$ & $59.2$ & $253.8$ & $193.8$ & $141.2$ & $44.1$ & \cellcolor{\myCellColor} $\mathbf{35.8}$ \\ 
qpG11 & $142.31$ & $7.3$ & $***$ & $4.79$ & \cellcolor{\myCellColor} $\mathbf{4.62}$ & $4.7$ & $305.1$ & $20.0$ & $12.3$ & \cellcolor{\myCellColor} $\mathbf{10.9}$ & $15.6$ & $13.2$ \\ 
qpG51 & $450.74$ & $186.49$ & $***$ & \cellcolor{\myCellColor} $\mathbf{89.81}$ & $150.86$ & $132.68$ & $523.6$ & $247.2$ & $246.9$ & $134.1$ & $73.9$ & \cellcolor{\myCellColor} $\mathbf{65.2}$ \\ 
thetaG11 & $332.06$ & $9.43$ & $***$ & $9.43$ & $9.46$ & \cellcolor{\myCellColor} $\mathbf{6.81}$ & $477.5$ & $21.9$ & $18.9$ & \cellcolor{\myCellColor} $\mathbf{12.5}$ & $16.3$ & $14.5$ \\ 
thetaG51 & $1062.91$ & $110.64$ & $***$ & $107.02$ & \cellcolor{\myCellColor} $\mathbf{37.22}$ & $85.92$ & $252.8$ & $230.8$ & $***\tnote{$\dagger$}$ & $131.7$ & $66.5$ & \cellcolor{\myCellColor} $\mathbf{53.7}$ \\ 
 \midrule 
rs1184 & $***\tnote{m}$ & $1227.48$ & $***$ & $882.27$ & $632.96$ & \cellcolor{\myCellColor} $\mathbf{569.29}$ & $***\tnote{m}$ & $4192.8$ & $3495.8$ & $3424.1$ & $2483.9$ & \cellcolor{\myCellColor} $\mathbf{2301.3}$ \\ 
rs1555 & $***\tnote{$\dagger$}$ & $79.83$ & $***$ & \cellcolor{\myCellColor} $\mathbf{65.93}$ & $80.72$ & $83.84$ & $***\tnote{$\dagger$}$ & $316.7$ & $242.7$ & $268.8$ & $160.7$ & \cellcolor{\myCellColor} $\mathbf{132.1}$ \\ 
rs1907 & $***\tnote{$\dagger$}$ & $233.86$ & $***$ & $197.99$ & $178.79$ & \cellcolor{\myCellColor} $\mathbf{166.23}$ & $***\tnote{$\dagger$}$ & $483.8$ & $490.3$ & $455.3$ & $382.7$ & \cellcolor{\myCellColor} $\mathbf{352.1}$ \\ 
rs200 & $640.0$ & $31.33$ & $***$ & $21.09$ & $24.78$ & \cellcolor{\myCellColor} $\mathbf{19.29}$ & $3366.2$ & $121.6$ & $93.0$ & $82.9$ & $93.2$ & \cellcolor{\myCellColor} $\mathbf{71.7}$ \\ 
rs228 & $206.2$ & $40.88$ & $***$ & $27.79$ & $25.29$ & \cellcolor{\myCellColor} $\mathbf{18.6}$ & $1220.2$ & $116.0$ & $59.4$ & $76.2$ & $67.1$ & \cellcolor{\myCellColor} $\mathbf{50.9}$ \\ 
rs35 & $296.52$ & $196.93$ & $***$ & $146.56$ & $88.86$ & \cellcolor{\myCellColor} $\mathbf{71.25}$ & $1269.8$ & $548.1$ & $358.2$ & $404.6$ & $272.5$ & \cellcolor{\myCellColor} $\mathbf{223.2}$ \\ 
rs365 & $***\tnote{$\dagger$}$ & $159.75$ & $***$ & $127.77$ & $110.48$ & \cellcolor{\myCellColor} $\mathbf{92.5}$ & $***\tnote{$\dagger$}$ & $433.1$ & $364.6$ & $351.1$ & $289.9$ & \cellcolor{\myCellColor} $\mathbf{262.0}$ \\ 
rs828 & $603.55$ & $29.86$ & $***$ & $19.24$ & $23.25$ & \cellcolor{\myCellColor} $\mathbf{17.81}$ & $3716.7$ & $113.2$ & $80.0$ & $71.1$ & $87.5$ & \cellcolor{\myCellColor} $\mathbf{64.2}$ \\ 

    \midrule
     problem & \multicolumn{6}{c}{Iterations} & \multicolumn{6}{c}{Number of cliques / Maximum clique size}\\
         \cmidrule(lr){2-7}\cmidrule(lr){8-13}
      & NoDe & NoMer & SpCo & ParCh & CG1 & CG2 &  NoDe & NoMer & SpCo & ParCh & CG1 & CG2 \\
     \midrule
    maxG11 & $225$ & $225$ & $500$ & $275$ & $200$ & $225$ & $1/800$ & $598/24$ & $13/80$ & $198/32$ & $473/28$ & $407/36$ \\ 
maxG32 & $300$ & $300$ & $425$ & $250$ & $225$ & $375$ & $1/2000$ & $1498/76$ & $21/210$ & $481/76$ & $1164/92$ & $664/102$ \\ 
maxG51 & $150$ & $100$ & $75$ & $50$ & $75$ & $125$ & $1/1000$ & $674/326$ & $181/322$ & $163/326$ & $448/362$ & $313/395$ \\ 
mcp500-1 & $350$ & $150$ & $150$ & $125$ & $125$ & $100$ & $1/500$ & $457/39$ & $451/44$ & $105/43$ & $437/54$ & $361/65$ \\ 
mcp500-2 & $225$ & $225$ & $200$ & $250$ & $150$ & $175$ & $1/500$ & $363/138$ & $144/138$ & $96/140$ & $316/156$ & $226/171$ \\ 
mcp500-3 & $200$ & $250$ & $250$ & $300$ & $200$ & $200$ & $1/500$ & $259/242$ & $101/242$ & $71/242$ & $211/263$ & $135/285$ \\ 
mcp500-4 & $225$ & $225$ & $375$ & $75$ & $100$ & $175$ & $1/500$ & $161/340$ & $63/346$ & $46/341$ & $105/368$ & $87/393$ \\ 
qpG11 & $400$ & $325$ & $525$ & $375$ & $250$ & $300$ & $1/1600$ & $1398/24$ & $813/80$ & $287/32$ & $1273/28$ & $1207/36$ \\ 
qpG51 & $750$ & $600$ & $800$ & $550$ & $1800$ & $1825$ & $1/2000$ & $1674/326$ & $1182/304$ & $275/326$ & $1448/362$ & $1313/395$ \\ 
thetaG11 & $675$ & $375$ & $2275$ & $650$ & $500$ & $400$ & $1/801$ & $598/25$ & $13/81$ & $198/33$ & $494/29$ & $423/41$ \\ 
thetaG51 & $3825$ & $325$ & $***\tnote{$\dagger$}$ & $575$ & $375$ & $1250$ & $1/1001$ & $676/324$ & $150/323$ & $157/324$ & $424/358$ & $267/396$ \\ 
 \midrule 
rs1184 & $***\tnote{m}$ & $225$ & $200$ & $200$ & $200$ & $200$ & $1/14822$ & $2236/500$ & $78/1330$ & $1043/500$ & $664/608$ & $258/632$ \\ 
rs1555 & $***\tnote{$\dagger$}$ & $150$ & $150$ & $150$ & $150$ & $175$ & $1/7479$ & $6891/184$ & $3350/187$ & $2556/184$ & $5529/236$ & $4858/276$ \\ 
rs1907 & $***\tnote{$\dagger$}$ & $200$ & $200$ & $175$ & $175$ & $200$ & $1/5357$ & $577/261$ & $47/585$ & $419/261$ & $441/324$ & $219/405$ \\ 
rs200 & $175$ & $125$ & $125$ & $125$ & $125$ & $125$ & $1/3025$ & $1632/95$ & $94/216$ & $444/95$ & $1123/112$ & $583/119$ \\ 
rs228 & $150$ & $125$ & $125$ & $125$ & $125$ & $125$ & $1/1919$ & $790/88$ & $48/180$ & $255/88$ & $369/95$ & $129/127$ \\ 
rs35 & $200$ & $175$ & $200$ & $175$ & $150$ & $150$ & $1/2003$ & $589/343$ & $53/735$ & $189/343$ & $214/457$ & $106/520$ \\ 
rs365 & $***\tnote{$\dagger$}$ & $175$ & $175$ & $175$ & $175$ & $175$ & $1/4704$ & $1230/296$ & $110/350$ & $539/296$ & $688/349$ & $346/474$ \\ 
rs828 & $150$ & $125$ & $150$ & $125$ & $125$ & $125$ & $1/3169$ & $1875/86$ & $112/174$ & $501/86$ & $1378/102$ & $708/126$ \\ 

   \end{tabular}}
 \begin{tablenotes}
 \scriptsize{
    \item[$\dagger$] time limit reached; \item[m] out of memory error; \item[1] no decomposition; \item[2] no merging; \item[3] \pkg{SparseCoLO} merging; \\ \item[4] parent-child merging;
    \item[5] clique graph with nominal edge weighting~\eqref{eq:edge_weighting}; \\ \item[6] clique graph with estimated edge weighting~\eqref{eq:estimated_edge_weighting}};
 \end{tablenotes}
 \end{threeparttable}
\end{center}
\end{table}

Our clique graph-based methods lead to a reduction in overall solver time.
The method with estimated edge weighting function \texttt{CG2} achieves the lowest average projection times for the majority of problems.
In four cases \texttt{ParCh} has a narrow advantage.
The geometric mean of the ratios of projection time of \texttt{CG2} compared to the best non-graph method is $\gmeanRatios$, with a minimum ratio of $\minRatio$ for problem \mbox{\texttt{\minRatioName}}.
There does not seem to be a clear pattern that relates the projection time to the number of cliques or the maximum clique size of the decomposition.
This is expected as the optimal merging strategy depends on the properties of the initial decomposition such as the overlap between the cliques.
The merging strategies \texttt{ParCh}, \texttt{CG1} and \texttt{CG2} generally result in similar maximum clique sizes compared to \texttt{SparseCoLO}, with \texttt{CG1} being the most conservative in the number of merges.

Table~\ref{tb:sdplib_solver_comparison} shows the benchmark results of \pkg{COSMO} with merging strategy \texttt{CG2}, \pkg{MOSEK}, and \pkg{SCS}.
The decomposition helps \pkg{COSMO} to solve most problems faster than \pkg{MOSEK} and \pkg{SCS}.
This is even more significant for the larger problems that were generated from the \pkg{SuiteSparse} Matrix Collection.
The decomposition does not seem to provide a major benefit for the slightly denser problems \texttt{mcp500-3} and \texttt{mcp500-4}.
Furthermore, \pkg{COSMO} seems to converge slowly for \texttt{qpG51} and \texttt{thetaG51}.
Similar observations for \texttt{mcp500-3}, \texttt{mcp500-4} and \texttt{thetaG51} have been made by the authors in~\cite{Andersen_2010a}.
Finally, many of the larger problems were not solvable within the time limit or caused out-of-memory problems if no decomposition was used in \pkg{MOSEK} and \pkg{SCS}.

\begin{table}[htb]
\sisetup{round-mode=places
        ,round-precision=2
        ,scientific-notation=true}
\begin{center}
\caption{Benchmark results for non-chordal sparse SDPs from SDPLib and SDPs generated with sparsity patterns from the \pkg{SuiteSparse} Matrix Collection.}
\label{tb:sdplib_solver_comparison}
\begin{threeparttable}
\resizebox{\textwidth}{!}{%
   \begin{tabular}{l l l l l l l l l l}
     \toprule
      & \multicolumn{3}{c}{Solve time (\si{s}) } & \multicolumn{3}{c}{Iterations} & \multicolumn{3}{c}{Max error\tnote{2}}\\
         \cmidrule(lr){2-4} \cmidrule(lr){5-7}\cmidrule(lr){8-10}
      problem & COSMO\tnote{1} & MOSEK & SCS &  COSMO\tnote{1} &  Mosek & SCS & COSMO\tnote{1} &  MOSEK & SCS\\
     \midrule
    maxG11 & \cellcolor{\myCellColor} $\mathbf{1.47}$ & $4.45$ & $131.8$ & 225 & 5 & 1220 & $ \num{0.0007842310829704409}$ & $ \num{0.0019802016970309646}$ & $ \num{0.0008740286728388414}$ \\ 
maxG32 & \cellcolor{\myCellColor} $\mathbf{6.25}$ & $50.84$ & $840.79$ & 375 & 5 & 1220 & $ \num{0.0009740208513682289}$ & $ \num{0.0027593976567979518}$ & $ \num{0.0006867004415874414}$ \\ 
maxG51 & \cellcolor{\myCellColor} $\mathbf{8.09}$ & $9.92$ & $36.56$ & 125 & 8 & 180 & $ \num{0.0023052421160915087}$ & $ \num{0.0002827382989427568}$ & $ \num{0.0008849426292532904}$ \\ 
mcp500-1 & \cellcolor{\myCellColor} $\mathbf{0.24}$ & $1.7$ & $29.28$ & 100 & 7 & 580 & $ \num{0.001019858373429112}$ & $ \num{0.001515028746100033}$ & $ \num{0.0006314096833755729}$ \\ 
mcp500-2 & \cellcolor{\myCellColor} $\mathbf{1.68}$ & $1.75$ & $17.36$ & 175 & 7 & 380 & $ \num{0.0006820071470928921}$ & $ \num{0.0007366481208091927}$ & $ \num{0.0003296179431987297}$ \\ 
mcp500-3 & $4.41$ & \cellcolor{\myCellColor} $\mathbf{1.68}$ & $8.36$ & 200 & 6 & 180 & $ \num{0.002226529296118806}$ & $ \num{0.0004179647134477934}$ & $ \num{0.0006862460288040091}$ \\ 
mcp500-4 & $8.2$ & \cellcolor{\myCellColor} $\mathbf{1.76}$ & $7.4$ & 175 & 7 & 160 & $ \num{0.00164832950036284}$ & $ \num{0.0002789769851946342}$ & $ \num{0.00034256415775601497}$ \\ 
qpG11 & \cellcolor{\myCellColor} $\mathbf{2.36}$ & $26.23$ & $734.7$ & 300 & 7 & 1820 & $ \num{0.00045656868412013237}$ & $ \num{0.0012781214356479534}$ & $ \num{0.0009766872920058663}$ \\ 
qpG51 & $121.6$ & \cellcolor{\myCellColor} $\mathbf{96.42}$ & $527.55$ & 1825 & 14 & 800 & $ \num{0.007558680070501833}$ & $ \num{0.0004802343556962799}$ & $ \num{0.0009827368649027318}$ \\ 
thetaG11 & \cellcolor{\myCellColor} $\mathbf{2.32}$ & $8.53$ & $142.53$ & 400 & 9 & 1380 & $ \num{0.0014337318382577717}$ & $ \num{0.0011103506672344022}$ & $ \num{8.310518951297808e-5}$ \\ 
thetaG51 & $71.21$ & \cellcolor{\myCellColor} $\mathbf{50.08}$ & $967.43$ & 1250 & 11 & 3240 & $ \num{0.13847976687437644}$ & $ \num{0.00040294361329714794}$ & $ \num{0.0009938423189167704}$ \\ 
\midrule 
rs1184 & \cellcolor{\myCellColor} $\mathbf{224.86}$ & $***\tnote{m}$ & $***\tnote{m}$ & 200 & *** & *** & $ \num{0.00056524162380964}$ & $ ***$ & $ ***$ \\ 
rs1555 & \cellcolor{\myCellColor} $\mathbf{66.6}$ & $***\tnote{$\dagger$}$ & $***\tnote{m}$ & 175 & *** & *** & $ \num{0.0005317204817832126}$ & $ ***$ & $ ***$ \\ 
rs1907 & \cellcolor{\myCellColor} $\mathbf{104.61}$ & $***\tnote{$\dagger$}$ & $***\tnote{m}$ & 200 & *** & *** & $ \num{0.00026645009322361464}$ & $ ***$ & $ ***$ \\ 
rs200 & \cellcolor{\myCellColor} $\mathbf{12.47}$ & $752.27$ & $***\tnote{$\dagger$}$ & 125 & 11 & *** & $ \num{0.00018692305793077034}$ & $ \num{0.00061060678892031}$ & $ ***$ \\ 
rs228 & \cellcolor{\myCellColor} $\mathbf{12.86}$ & $395.24$ & $982.5$ & 125 & 11 & 1620 & $ \num{0.00021461744055341749}$ & $ \num{0.0006342729996048659}$ & $ \num{0.0008268745300744567}$ \\ 
rs35 & \cellcolor{\myCellColor} $\mathbf{54.88}$ & $919.19$ & $***\tnote{$\dagger$}$ & 150 & 12 & *** & $ \num{0.00024752587359667803}$ & $ \num{0.00032619870648789784}$ & $ ***$ \\ 
rs365 & \cellcolor{\myCellColor} $\mathbf{62.65}$ & $***\tnote{$\dagger$}$ & $***\tnote{$\dagger$}$ & 175 & *** & *** & $ \num{0.0003079972836244303}$ & $ ***$ & $ ***$ \\ 
rs828 & \cellcolor{\myCellColor} $\mathbf{10.84}$ & $825.03$ & $***\tnote{$\dagger$}$ & 125 & 11 & *** & $ \num{0.00019847358220930543}$ & $ \num{0.0008480747078987047}$ & $ ***$ \\ 

   \bottomrule
   \end{tabular}}
 \begin{tablenotes}
 \scriptsize{%
    \item[1] with chordal decomposition and clique merging strategy \texttt{CG2}; \item[2] $\max\{\epsilon_1, \epsilon_2, \epsilon_3\}$;
    \item[$\dagger$] time limit reached;\\
     \item[m] out of memory error;
    }%
 \end{tablenotes}
 \end{threeparttable}
\end{center}
\end{table}


\section{Conclusions}
\label{sec:conclusion}
This paper describes the first-order solver \pkg{COSMO} and the ADMM algorithm on which it is based.
The solver combines direct support of quadratic objectives, infeasibility detection, custom constraints, chordal decomposition of PSD constraints and automatic clique merging.
The performance of the solver is illustrated on a number of benchmark problems that challenge different aspects of modern solvers.

The implementation in the Julia language facilitates rapid development and testing of ideas and allows users to customize the solver for their applications.
It further allows the abstraction of precision and array types which we are planning to use to allow \pkg{COSMO} to run on GPUs.
Further performance gains are likely to be achieved by exploring acceleration methods to speed up convergence to higher accuracies and reduce the dependency on problem scaling.

\newpage
\appendix

\section{Appendix}

\subsection{Benchmark Results}
\FloatBarrier

\begin{table}[htb]
\sisetup{round-mode=places
        ,round-precision=2
        ,scientific-notation=true}
\begin{center}
\caption{Benchmark results for nearest correlation matrix problems.}
\label{tb:correlation}
\begin{threeparttable}
\resizebox{\textwidth}{!}{%
   \begin{tabular}{l l l l l l l l l l}
     \toprule
      & \multicolumn{3}{c}{Solve time (\si{s}) } & \multicolumn{3}{c}{Iterations} & \multicolumn{3}{c}{Max error\tnote{2}}\\
         \cmidrule(lr){2-4} \cmidrule(lr){5-7}\cmidrule(lr){8-10}
      n & COSMO & MOSEK & SCS &  COSMO &  MOSEK & SCS & COSMO &  MOSEK & SCS\\
     \midrule
    $50$  & $0.5$ & $1.89$ & \cellcolor{\myCellColor} $\mathbf{0.07}$ & 25 & 4 & 20 & $ \num{2.0477763720282582e-5}$ & $ \num{0.0013520741790252007}$ & $ \num{6.183471303045888e-7}$ \\ 
$100$  & $0.45$ & $35.68$ & \cellcolor{\myCellColor} $\mathbf{0.18}$ & 25 & 4 & 20 & $ \num{2.024643461273149e-5}$ & $ \num{0.0062584935298491615}$ & $ \num{4.244552375162742e-6}$ \\ 
$150$  & \cellcolor{\myCellColor} $\mathbf{0.39}$ & $244.35$ & $0.46$ & 25 & 4 & 20 & $ \num{5.1477147785084844e-5}$ & $ \num{0.0086072456008019}$ & $ \num{9.001628710297762e-6}$ \\ 
$200$  & $0.94$ & $1032.25$ & \cellcolor{\myCellColor} $\mathbf{0.74}$ & 25 & 4 & 20 & $ \num{7.586271996820478e-5}$ & $ \num{0.003641904488103727}$ & $ \num{8.412082256913692e-5}$ \\ 
$400$  & \cellcolor{\myCellColor} $\mathbf{1.14}$ & $***\tnote{$\dagger$}$ & $3.07$ & 25 & *** & 20 & $ \num{0.0002620943755379533}$ & $ ***$ & $ \num{8.521790746760832e-5}$ \\ 
$600$  & \cellcolor{\myCellColor} $\mathbf{3.05}$ & $***\tnote{$\dagger$}$ & $7.38$ & 25 & *** & 20 & $ \num{0.00015342773517135184}$ & $ ***$ & $ \num{0.00010114042713642336}$ \\ 
$800$  & \cellcolor{\myCellColor} $\mathbf{5.08}$ & $***\tnote{$\dagger$}$ & $13.85$ & 25 & *** & 20 & $ \num{0.00040339879168537667}$ & $ ***$ & $ \num{0.00010736567516677556}$ \\ 
$1000$  & \cellcolor{\myCellColor} $\mathbf{8.43}$ & $***\tnote{$\dagger$}$ & $21.94$ & 25 & *** & 20 & $ \num{0.0008107455418383952}$ & $ ***$ & $ \num{0.00014040012268331096}$ \\ 
$1200$  & \cellcolor{\myCellColor} $\mathbf{12.1}$ & $***\tnote{$\dagger$}$ & $34.58$ & 25 & *** & 20 & $ \num{0.001023757730316243}$ & $ ***$ & $ \num{0.00017292648579278475}$ \\ 
$1400$  & \cellcolor{\myCellColor} $\mathbf{17.19}$ & $***\tnote{$\dagger$}$ & $49.23$ & 25 & *** & 20 & $ \num{0.0010223037071205527}$ & $ ***$ & $ \num{0.0001929429694929055}$ \\ 
$1600$  & \cellcolor{\myCellColor} $\mathbf{22.69}$ & $***\tnote{$\dagger$}$ & $61.07$ & 25 & *** & 20 & $ \num{0.0008002665862960172}$ & $ ***$ & $ \num{0.00022140571363724054}$ \\ 
$1800$  & \cellcolor{\myCellColor} $\mathbf{29.41}$ & $***\tnote{$\dagger$}$ & $74.0$ & 25 & *** & 20 & $ \num{0.0005144474574880278}$ & $ ***$ & $ \num{0.0002434220709377156}$ \\ 
$2000$  & \cellcolor{\myCellColor} $\mathbf{37.74}$ & $***\tnote{$\dagger$}$ & $101.18$ & 25 & *** & 20 & $ \num{0.0003312357476843221}$ & $ ***$ & $ \num{0.00025673549576654914}$ \\ 

   \bottomrule
   \end{tabular}}
  \begin{tablenotes}
  \scriptsize{%
     \item[2] $\max\{\epsilon_1, \epsilon_2, \epsilon_3\}$;
     \item[$\dagger$] time limit reached
     }%
  \end{tablenotes}
 \end{threeparttable}
\end{center}
\end{table}

\begin{table}[htb]
\begin{center}
\sisetup{round-mode=places
        ,round-precision=2
        ,scientific-notation=true}
\caption{Benchmark results for block arrow sparse SDPs with varying number of blocks $N_b$ and block size $d$.}\label{tb:blockarrow}
\begin{threeparttable}
\resizebox{\textwidth}{!}{%
   \begin{tabular}{l l l l l l l l l l l l l}
     \toprule
      & \multicolumn{4}{c}{Solve time (\si{s}) } & \multicolumn{4}{c}{Iterations} & \multicolumn{4}{c}{Max error\tnote{2}}\\
         \cmidrule(lr){2-5} \cmidrule(lr){6-9}\cmidrule(lr){10-13}
      $d$ & COSMO & COSMO\tnote{cd} &  MOSEK & SCS & COSMO & COSMO\tnote{cd} & MOSEK & SCS & COSMO & COSMO\tnote{cd} &  MOSEK & SCS\\
     \midrule
    $10$  & $12.71$ & \cellcolor{\myCellColor} $\mathbf{1.22}$ & $16.5$ & $56.99$ & 125 & 175 & 9 & 1340 & $ \num{0.00015007208583617837}$ & $ \num{0.0005599589215771029}$ & $ \num{0.0011644230211438464}$ & $ \num{0.0009007495080682709}$ \\ 
$12$  & $23.95$ & \cellcolor{\myCellColor} $\mathbf{2.3}$ & $28.72$ & $68.72$ & 125 & 225 & 10 & 1120 & $ \num{0.00022915625553989653}$ & $ \num{0.0007520071703709815}$ & $ \num{3.6530307821706646e-5}$ & $ \num{0.0005379673419291557}$ \\ 
$14$  & $27.4$ & \cellcolor{\myCellColor} $\mathbf{2.39}$ & $41.66$ & $147.72$ & 150 & 175 & 10 & 1780 & $ \num{7.915695395411079e-5}$ & $ \num{0.0008601620115847063}$ & $ \num{0.000684317980750578}$ & $ \num{0.0005461080386119225}$ \\ 
$16$  & $36.4$ & \cellcolor{\myCellColor} $\mathbf{3.6}$ & $58.04$ & $254.93$ & 150 & 275 & 10 & 2360 & $ \num{4.755751338248006e-5}$ & $ \num{0.0005340590506608241}$ & $ \num{0.00026373003717226066}$ & $ \num{0.0006913103908292647}$ \\ 
$18$  & $42.96$ & \cellcolor{\myCellColor} $\mathbf{4.42}$ & $85.81$ & $472.77$ & 150 & 275 & 11 & 3440 & $ \num{0.0002662293488024008}$ & $ \num{0.000576826639568577}$ & $ \num{0.00025482498041929497}$ & $ \num{0.000325376397877594}$ \\ 
$20$  & $41.13$ & \cellcolor{\myCellColor} $\mathbf{5.08}$ & $112.56$ & $251.39$ & 125 & 275 & 11 & 1460 & $ \num{0.00026081829121172055}$ & $ \num{0.0005347945419812834}$ & $ \num{6.328217577546178e-5}$ & $ \num{0.000726192018424364}$ \\ 
$22$  & $50.7$ & \cellcolor{\myCellColor} $\mathbf{11.47}$ & $142.14$ & $520.81$ & 125 & 700 & 11 & 2460 & $ \num{0.0002867082609768577}$ & $ \num{0.0005845225291909384}$ & $ \num{0.0003310618475841763}$ & $ \num{0.0004857063616828414}$ \\ 
$24$  & $125.12$ & \cellcolor{\myCellColor} $\mathbf{10.52}$ & $178.75$ & $445.37$ & 275 & 500 & 11 & 1780 & $ \num{7.455716137394613e-5}$ & $ \num{0.0007232705962837329}$ & $ \num{0.00035676724430774604}$ & $ \num{0.000509563519333417}$ \\ 
$26$  & $88.74$ & \cellcolor{\myCellColor} $\mathbf{6.84}$ & $220.26$ & $623.46$ & 150 & 200 & 11 & 2060 & $ \num{8.169110186283036e-5}$ & $ \num{0.00032748159042416086}$ & $ \num{2.8875185195929123e-6}$ & $ \num{0.00013852797612205354}$ \\ 
$28$  & $98.83$ & \cellcolor{\myCellColor} $\mathbf{7.47}$ & $263.52$ & $1626.82$ & 150 & 200 & 11 & 4660 & $ \num{3.7339497199617757e-5}$ & $ \num{0.0008925344498322185}$ & $ \num{2.6174549682816456e-5}$ & $ \num{0.00045919955870531805}$ \\ 

    \midrule
    $N_b$ & COSMO & COSMO\tnote{cd} &  MOSEK & SCS & COSMO & COSMO\tnote{cd} & MOSEK & SCS & COSMO & COSMO\tnote{cd} &  MOSEK & SCS\\
    \midrule
     $50$  & $17.31$ & \cellcolor{\myCellColor} $\mathbf{3.73}$ & $28.69$ & $115.36$ & 150 & 300 & 10 & 2540 & $ \num{0.0001326866826418272}$ & $ \num{0.0006283058996429108}$ & $ \num{0.00022114584945078286}$ & $ \num{0.0009121165420083106}$ \\ 
$60$  & $21.66$ & \cellcolor{\myCellColor} $\mathbf{6.4}$ & $46.38$ & $217.75$ & 150 & 475 & 11 & 3380 & $ \num{0.00014235168846067205}$ & $ \num{0.0008741453569003996}$ & $ \num{7.271752765152321e-5}$ & $ \num{0.0008814324823583168}$ \\ 
$70$  & $34.51$ & \cellcolor{\myCellColor} $\mathbf{9.89}$ & $64.8$ & $177.24$ & 175 & 700 & 11 & 2060 & $ \num{4.121757838328223e-5}$ & $ \num{0.0006136015839921192}$ & $ \num{1.7164282938544684e-5}$ & $ \num{0.0008992284546667715}$ \\ 
$80$  & $37.22$ & \cellcolor{\myCellColor} $\mathbf{5.78}$ & $85.29$ & $175.38$ & 175 & 275 & 11 & 1540 & $ \num{5.1520246418980796e-5}$ & $ \num{0.0007132984060971376}$ & $ \num{0.00020795501333830744}$ & $ \num{0.0004566101336938731}$ \\ 
$90$  & $63.32$ & \cellcolor{\myCellColor} $\mathbf{12.85}$ & $113.24$ & $234.56$ & 225 & 550 & 11 & 1580 & $ \num{4.800650126711851e-5}$ & $ \num{0.0011211326450515217}$ & $ \num{4.788398188007093e-5}$ & $ \num{0.0006403440946790319}$ \\ 
$100$  & $67.3$ & \cellcolor{\myCellColor} $\mathbf{25.16}$ & $140.81$ & $328.2$ & 200 & 1100 & 11 & 1800 & $ \num{6.383572514351841e-5}$ & $ \num{0.000953976307407238}$ & $ \num{3.277052034139081e-5}$ & $ \num{0.0008593130862501243}$ \\ 
$110$  & $117.7$ & \cellcolor{\myCellColor} $\mathbf{39.86}$ & $172.99$ & $761.22$ & 275 & 1650 & 11 & 3140 & $ \num{3.7775769869920564e-5}$ & $ \num{0.0006800982666091884}$ & $ \num{0.00013501200857849322}$ & $ \num{0.0003526527965693732}$ \\ 
$120$  & $81.31$ & \cellcolor{\myCellColor} $\mathbf{29.06}$ & $209.96$ & $1191.64$ & 150 & 1050 & 11 & 4620 & $ \num{0.0003591889356256073}$ & $ \num{0.0010699512965270312}$ & $ \num{0.0003660384318623842}$ & $ \num{0.00027095669266488805}$ \\ 
$130$  & $115.21$ & \cellcolor{\myCellColor} $\mathbf{26.61}$ & $269.06$ & $1094.53$ & 175 & 925 & 12 & 3580 & $ \num{8.482631598897952e-5}$ & $ \num{0.0013118412845728729}$ & $ \num{8.592512773412302e-5}$ & $ \num{0.0009669220513205416}$ \\ 
$140$  & $114.95$ & \cellcolor{\myCellColor} $\mathbf{26.87}$ & $290.23$ & $***\tnote{$\dagger$}$ & 150 & 800 & 11 & *** & $ \num{0.00011982337858651794}$ & $ \num{0.0019225218705366573}$ & $ \num{0.000498701306519267}$ & $ ***$ \\ 

   \bottomrule
   \end{tabular}}

 \begin{tablenotes}
 \scriptsize{%
    \item[cd] with chordal decomposition; \item[2] $\max\{\epsilon_1, \epsilon_2, \epsilon_3\}$; \item[$\dagger$] time limit reached
    }%
 \end{tablenotes}
 \end{threeparttable}
\end{center}
\end{table}

\FloatBarrier
\normalsize
\section{Appendix}
\subsection{Projection onto $\mathcal{C}_{\sum}$}
\label{app_project}
The projection of a symmetric matrix $S \in \Sym[n]$ onto the set of matrices where the sum of rows and columns each equal to one can be written as the following optimisation problem:
\begin{equation}
  \begin{array}{ll}
    \mbox{minimize}   & \textstyle{\frac{1}{2}} \norm{s_p - s}_F^2\\
    \mbox{subject to} & A s_p = \mathbf{1}_{2n-1}
  \end{array}
  \label{eq:proj_dsm}
\end{equation}
where $s = \text{vec}(S)$ is the vectorized matrix, $s_p$ is the projected vector and $A$ is given by:
\begin{equation*}
A = \begin{bmatrix} A_r & A_c \end{bmatrix}^\top
\end{equation*}
where $A_r$ is used to constrain the rows of $S$ and $A_c$ is used to constrain the columns of $S$:
\begin{equation*}
A_r = \mathbf{1}_n^\top \otimes I_n \quad A_c = \begin{bmatrix} I_{n-1} \otimes \mathbf{1}_n^\top & \mathbf{0}_{n-1 \times n} \end{bmatrix}.
\end{equation*}
Notice that $A_c$ is a $(n-1) \times n$ matrix because the redundant constraint on the last column was removed.
By writing down the KKT optimality conditions for~\eqref{eq:proj_dsm} we obtain the following equations:
\begin{align}
A s_p  &= \mathbf{1}_{2n-1}, \label{eq:dsm_kkt1}\\
 s_p - s + A^\top \eta &= 0,\label{eq:dsm_kkt2}
\end{align}
with dual variable $\eta$. Eliminating $s_p$ and solving for $\eta$ yields
\begin{equation}
\label{eq:eta}
  \eta = (A A^\top)^{-1} (A s - \mathbf{1}_{2n-1}).
\end{equation}
The projected vector $s_p$ can be recovered from~\eqref{eq:dsm_kkt2}:
\begin{equation}
\label{eq:final_proj}
  s_p = s - A^\top \eta.
\end{equation}
It turns out that the inverse of $AA^\top$ in~\eqref{eq:eta} can be efficiently computed without a factorisation. To see this, form $AA^\top$ and write~\eqref{eq:eta} as:
\begin{equation}
\label{eq:eta_matrix}
  \begin{bmatrix}
  n I_n & \mathbf{1}_n \mathbf{1}_{n-1}^\top \\ \mathbf{1}_{n-1} \mathbf{1}_{n}^\top & n I_{n-1}
  \end{bmatrix}
  \begin{bmatrix}
    \eta_1 \\ \eta_2
  \end{bmatrix}
  = \begin{bmatrix}
  r_1 \\ r_2
  \end{bmatrix} \, \text{with }  \begin{bmatrix}
  r_1 \\ r_2
  \end{bmatrix} = A s - \mathbf{1}_{2n-1},
\end{equation}
where $\eta$ and the right-hand side where partitioned in such a way that $\eta_1, r_1 \in \Re^n$ and $\eta_2, r_2 \in \Re^{n-1}$.
Eliminating $\eta_1$ from~\eqref{eq:eta_matrix}, yields:
\begin{equation}
  \left(I_{n-1} - \textstyle{\frac{1}{n}} \mathbf{1}_{n-1} \mathbf{1}_{n-1}^\top  \right) \eta_2 =\textstyle{\frac{1}{n}} \left( r_2 - \textstyle{\frac{1}{n}} \mathbf{1}_{n-1} \mathbf{1}_{n}^\top  r_1 \right)
\end{equation}
The vector $\eta_2$ can be computed by applying the Sherman-Morrison formula to compute the inverse of the matrix on the left:
\begin{equation}
  \eta_2 = \textstyle{\frac{1}{n}} \left( I_{n-1} + \mathbf{1}_{n-1} \mathbf{1}_{n-1}^\top \right) \left( r_2 - \textstyle{\frac{1}{n}} \mathbf{1}_{n-1} \mathbf{1}_{n}^\top r_1 \right).
\end{equation}
The vector $\eta_1$ is then computed by substituting into the upper equation in~\eqref{eq:eta_matrix}:
\begin{equation}
  \eta_1 = \textstyle{\frac{1}{n}} \left(r_1 - \mathbf{1}_{n} \mathbf{1}_{n-1}^\top \eta_2 \right).
\end{equation}
Having computed the elements of the dual variable $\eta$, the projected vector $s_p$ is obtained by solving~\eqref{eq:final_proj}.

\bibliography{COSMO_bibliography.bib}

\end{document}